\documentclass[12pt,twoside,sort&compress ]{article}
\usepackage{mathrsfs,amsfonts,amsmath,amsbsy}\usepackage{fontenc}\usepackage{textcomp}
\usepackage{graphicx}
\usepackage{amssymb}
\usepackage{graphicx}

\def\le{\leq}

\def\bbb{\begin{eqnarray*}}

\def\eee{\end{eqnarray*}}

\begin{document}
\baselineskip=18pt
\begin{center}
\vspace{-0.6in}{\large \bf Continuous Dependence of the $n$-th Eigenvalue on \\[0.1in]
 Self-adjoint Discrete  Sturm-Liouville Problem}
\\ [0.2in]
 Hao Zhu,\ \ Yuming Shi $^{\dag}$  \\
\vspace{0.15in}  Department of Mathematics, Shandong University\\
Jinan, Shandong 250100, P. R. China\\

\footnote{$^\dag$ The corresponding author.}
\footnote{ \;\;Email addresses: haozhusdu@163.com(H. Zhu), ymshi@sdu.edu.cn(Y. Shi).}
\end{center}

{\bf Abstract.} This paper is concerned with continuous dependence of the $n$-th eigenvalue on self-adjoint discrete Sturm-Liouville
 problems.  The $n$-th eigenvalue is considered as a function in  the space of the problems.
A necessary and sufficient condition for all the eigenvalue functions to be continuous  and several
properties of the  eigenvalue functions in a set of the space of the problems are given.
They play an important role in the study of   continuous dependence of the $n$-th eigenvalue function on the problems.
Continuous dependence of  the $n$-th eigenvalue function  on the  equations and on the boundary conditions is studied separately. Consequently, the continuity and discontinuity of the $n$-th eigenvalue function are completely characterized in the whole space of the
 problems. Especially, asymptotic behaviors of the $n$-th eigenvalue function near each discontinuity point are given.
\medskip

\medskip

{2000 {\bf \it Mathematics Subject Classification}}: 39A12; 34B24;
39A70.

\bigskip

\noindent{\bf 1. Introduction}\medskip

A self-adjoint discrete Sturm-Liouville problem (briefly, SLP) considered in the present paper consists of a symmetric discrete Sturm-Liouville equation (briefly, SLE) \vspace{-0.2cm}
$$ -\nabla(f_{n}\Delta y_{n})+q_{n}y_{n}=\lambda w_{n}y_{n}, \;\;\;\; \ \  n\in[1,N],                                                                 \eqno(1.1)
\vspace{-0.2cm} $$
 and a self-adjoint boundary condition (briefly, BC)\vspace{-0.2cm}
$$
A \left (\begin{array} {cc}y_{0}\\
f_{0}\triangle y_{0}\end{array}   \right )
+B\left ( \begin{array} {cc}y_{N}\\
f_{N}\triangle y_{N}
\end{array}   \right )=0,                                                                                                                             \eqno(1.2)
\vspace{-0.2cm}
$$
where $N \geq 2$ is an integer, $\Delta $ and $\nabla$ are the forward and backward difference
operators, respectively, i.e., $\Delta y_n=y_{n+1}-y_n$ and $\nabla y_n=y_n-y_{n-1}$;
$f=\{f_n\}_{n=0}^{N}$, $q=\{q_n\}_{n=1}^{N}$ and $w=\{w_n\}_{n=1}^{N}$ are real-valued sequences such that
\vspace{-0.2cm}$$f_n\neq0 \;\;{\rm for}\;\; n\in [0,N],
\;\;w_n>0\;\;{\rm for}\;\; n\in [1,N];                                                                                                                 \eqno(1.3)
\vspace{-0.2cm}$$
$\lambda$ is the spectral parameter; the interval $[M,N]$ denotes the set of integers
 $\{M,M+1,\cdots,N\}$; $A$ and $B$ are
$2\times2$ complex matrices such that rank(A,B)=2, and satisfy the following self-adjoint boundary condition:
\vspace{-0.2cm} $$
 A\left (\begin{array} {cc} 0&1\\
-1&0\end{array}   \right )A^{*}=
B\left (\begin{array} {cc} 0&1\\
-1&0\end{array}   \right )B^{*},\vspace{-0.2cm}                                                                                                       \eqno(1.4)
$$
while $A^{*}$ denotes the complex conjugate transpose of $A$.

 Throughout this paper, by $\mathbb{C}$, $\mathbb{R}$, and $\mathbb{N}$ denote the sets of the complex, real, and natural numbers, respectively; and by $\bar{z}$ denotes the
complex conjugate of $z\in \mathbb{C}$. When a capital Latin letter stands for a
matrix, its entries  are denoted by the corresponding
lower case letter with two indices. For example, the entries of a
matrix $C$ are $c_{ij}$'s.

The SLP (1.1)--(1.2) has attracted increasing attention by many scholars, for its  applications in various areas of physical science [1, 9]. The $n$-th eigenvalue of the SLP always stands for an important value in the physical problems and it varies as the SLP varies.
A natural question is how it varies as the SLP varies. In the present paper,  we are interested in  continuous dependence of the $n$-th eigenvalue on the SLP.
We shall consider the $n$-th eigenvalue as a function in the space of the SLPs, and mainly study its continuity and discontinuity, and characterize
its asymptotic behaviors near each discontinuity point.

The continuous dependence of the $n$-th eigenvalue on self-adjoint continuous Sturm-Liouville problems
has been studied quite deeply and some elegant results have been obtained (cf., [6, 8, 10, 11, 14, 21]).
Now, we shall briefly recall some existing results of continuous dependence of the $n$-th eigenvalue on self-adjoint continuous SLPs. A self-adjoint continuous SLP consists of a differential SLE
$$ -(p(t)y')'+q(t)y=\lambda w(t)y, \;\;\ \  t\in(a,b),                       \eqno(1.5) \vspace{-0.2cm} $$
and a BC \vspace{-0.2cm}
$$
A \left (\begin{array} {cc}y(a)\\
(py')(a)\end{array}   \right )
+B\left ( \begin{array} {cc}y(b)\\
(py')(b)
\end{array}   \right )=0,                                                 \eqno(1.6)\vspace{-0.05cm}$$
where $-\infty<a<b<+\infty$; $1/p,q,w\in L((a,b),\mathbb{R})$ with $p,w>0$
almost everywhere  in $(a,b)$, while $L((a,b),\mathbb{R})$ denotes the space
of Lebesgue integrable
real functions in $(a,b)$; $A$ and $B$ are
$2\times2$ complex matrices such that rank(A,B)=2 and (1.4) holds.
It is well-known that the problem (1.5)--(1.6) has infinitely countable eigenvalues,
which are all real and can be
arranged  in the following non-decreasing order:
\vspace{-0.2cm}$$\lambda_0\le \lambda_1 \le \lambda_2 \le \cdots \le \lambda_n\le \cdots \eqno(1.7)\vspace{-0.2cm}$$
with $\lambda_n\to +\infty$ as $n\to +\infty$, counting repeatedly according to
their multiplicities. Using the variational method,  Courant and  Hilbert  showed that
 the $n$-th eigenvalue $\lambda_n$ is continuously dependent  on the problem
under the assumptions that the coefficient functions $p$, $q$ and $w$ in (1.5) are continuous functions in $(a,b)$
and the BC (1.6) is of a special class [6].
In 1997, Everitt, M\"{o}ller, and Zettl showed that  the $n$-th eigenvalue $\lambda_n$ does not depend  continuously on the BCs in general [8].
Later, Kong, Wu, and Zettl deeply studied this problem in 1999 [10]. They showed that
the $n$-th eigenvalue $\lambda_n$ depends continuously on the SLEs, found its discontinuity set
of the space of BCs
\vspace{-0.05cm}$$\begin{array} {cccc}\mathcal{J}^\mathbb{C}=\left\{[e^{i\theta}K|-I]:\,K\in SL(2,\mathbb{R}),k_{12}=0,\theta\in[0,\pi)\right\}\\[1.0ex]
\bigcup\left\{\left [\begin{array} {cccc}a_1&a_2&0&0\\
0&0&b_1&b_2\end{array}  \right ]:\,a_2b_2=0\right\},               \end{array}\eqno(1.8)\vspace{-0.2cm}$$
and gave its asymptotic behaviors near each discontinuity point.

The spectral theory of self-adjoint discrete Sturm-Liouville problems has attracted a great deal of interest
from many authors and some good results have been obtained (cf., [1, 3--5, 9, 12, 15--20, 22]).
In [17], the second author of the present paper with her coauthor Chen  studied some regular self-adjoint spectral problems for second-order vector difference equations, which include (1.1)--(1.2), and gave several spectral results, including the reality of the eigenvalues, the finiteness of the number of the eigenvalues,  and a formula for counting the number of the eigenvalues.
Based on these results, the problem (1.1)--(1.2) has $k$ eigenvalues (counting multiplicities), which are real and   can be arranged in the following non-decreasing order:
\vspace{-0.2cm}$$\lambda_0\le \lambda_1\le \lambda_2 \le \cdots\le\lambda_{k-1},\vspace{-0.2cm}$$
where $k$ can be determined (see Lemma 2.4). Note that the analytic and geometric multiplicities of an eigenvalue of an SLP (1.1)--(1.2) are the same [18, 22]. Recently, we studied some problems about dependence of the eigenvalues of (1.1)--(1.2) on the problems in [22]. We gave the topologies and geometric structures of the space of the SLPs (1.1)--(1.2), showed that each eigenvalue of a given SLP lives in one or two continuous eigenvalue branches, and  studied analyticity, differentiability and monotonicity of continuous eigenvalue branches.

It is evident that the $n$-th eigenvalue depends on the SLP (1.1)--(1.2)  and then can be regarded as a function in  the space of the SLPs. So, based on the work in [22], we shall mainly study continuous dependence of the $n$-th eigenvalue on the problem in the present paper. We shall give out its continuity and discontinuity sets in the space of
the SLPs (1.1)--(1.2), and characterize its asymptotic behaviors near each discontinuity point.

From Example 5.3 in [22], we have found that  the
index of eigenvalues in a continuous eigenvalue branch  may change as the
problem varies.  This may lead the discontinuity of the $n$-th eigenvalue function.
For convenience, we shall briefly recall this example.\medskip

\noindent{\bf Example 1.1} [22, Example 5.3]. Consider (1.1)--(1.2), where
\vspace{-0.2cm}$$
N=2, \; f_0 =f_1=f_2 =1, \; q_1 =q_2 =0, \; w_1 = w_2 =1,
\vspace{-0.2cm}$$
and   \vspace{-0.2cm}$$A=A(\alpha)=\left (\begin{array} {cccc}\cos\alpha&-\sin\alpha\\0&0\end{array}  \right ),\;
B=\left (\begin{array} {cccc}0&0\\0&-1\end{array}  \right ).\vspace{-0.05cm}$$
We have showed that the SLP with $\alpha=3\pi/4$ has exactly one eigenvalue $\lambda_0
=1$, and  each SLP with $\alpha \in [0,3\pi/4) \cup (3\pi/4,\pi)$ has exactly the following two eigenvalues:
\vspace{-0.2cm}$$
\lambda_0(\alpha) =\begin{cases} \lambda_-(\alpha) & \text{ if }
\alpha \in [0,3\pi/4),
     \\ \lambda_+(\alpha) & \text{ if } \alpha \in (3\pi/4,\pi),
     \end{cases}
\qquad \lambda_1(\alpha) =\begin{cases} \lambda_+(\alpha) & \text{
if } \alpha \in [0,3\pi/4),
     \\ \lambda_-(\alpha) & \text{ if } \alpha \in (3\pi/4,\pi),
     \end{cases}
\vspace{-0.2cm}$$
where
\vspace{-0.2cm}$$
\lambda_{\pm}(\alpha) ={3\cos\alpha+2\sin\alpha \pm
\sqrt{\cos^2\alpha+4\sin(2\alpha)+4} \over
  2(\cos\alpha+\sin\alpha) }.
$$
So, there are exactly  the following three continuous eigenvalue branches:
\vspace{-0.2cm}$$\begin{array}{ll}
\Lambda_1(\alpha) & =
    \lambda_0(\alpha) \text{ \ for } \alpha \in [0,3\pi/4),
\\
\Lambda_{2}(\alpha) & = \begin{cases} \lambda_1(\alpha) &
\text{ if } \alpha \in [0,3\pi/4),
    \\ 1 & \text{ if } \alpha =3\pi/4,
    \\ \lambda_0(\alpha) & \text{ if } \alpha \in (3\pi/4,\pi),
    \end{cases}
\\
\Lambda_3(\alpha) & = \lambda_1(\alpha) \text{ \ for } \alpha \in
(3\pi/4,\pi).
\end{array}\vspace{-0.2cm}$$
See Figure 1.1.
\vspace{-0.5cm}
\begin{center}
\includegraphics[width=87mm]{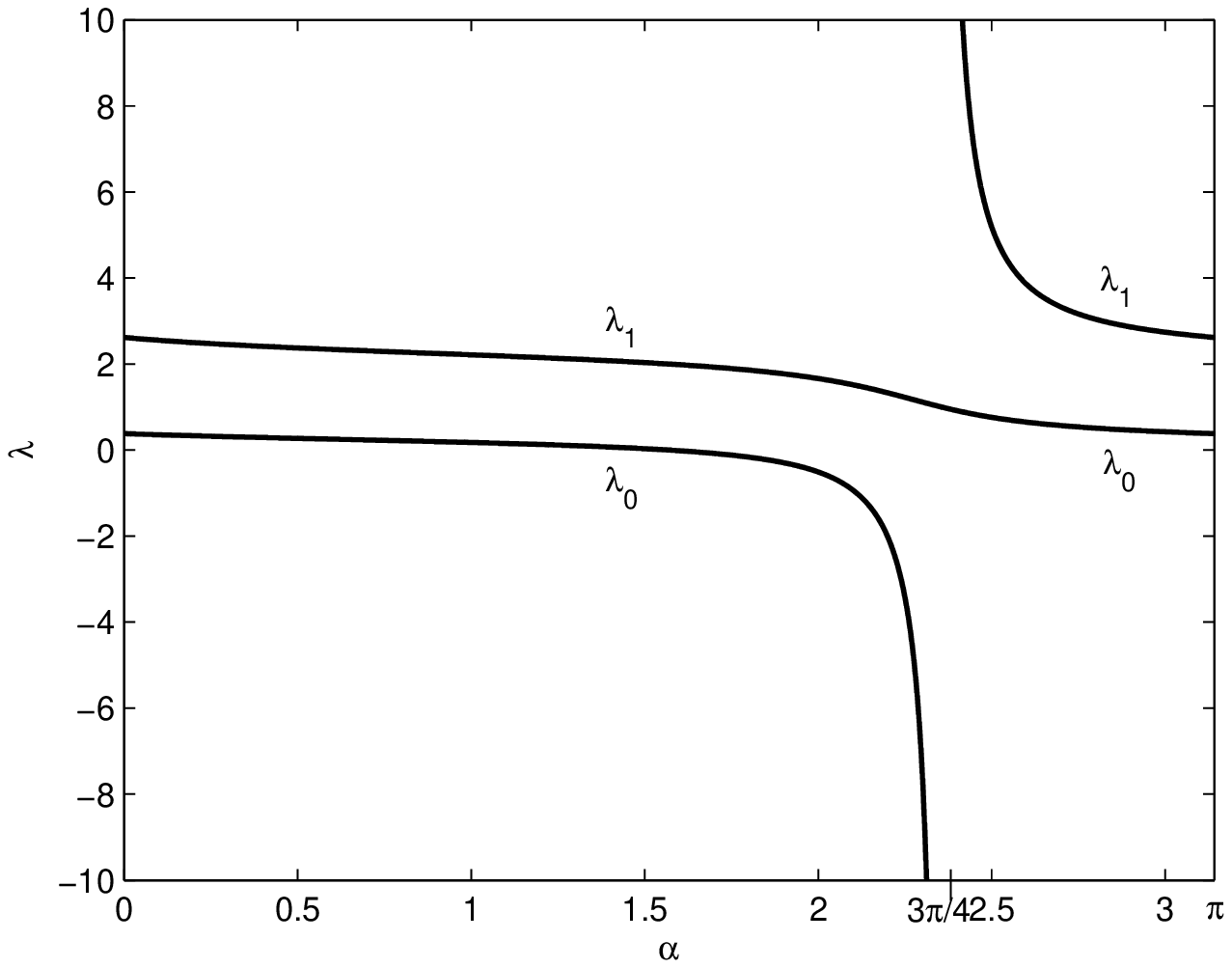}

{\small {\bf Figure 1.1.} }
\end{center}

As  functions in $\alpha\in[0,\pi)$, the eigenvalues $\lambda_0$ and $\lambda_1$ are not continuous at $\alpha=3\pi/4$, and have  the following asymptotic behaviors near $3\pi/4$:
\vspace{-0.2cm}$$\begin{array} {llll}\lim\limits_{\alpha\rightarrow 3\pi/4^-}\lambda_{0}(\alpha)=-\infty,\;\;
\lim\limits_{\alpha\rightarrow 3\pi/4^-}\lambda_{1}(\alpha)=\lambda_{0}(3\pi/4),\\[1.0ex]
\lim\limits_{\alpha\rightarrow 3\pi/4^+}\lambda_{0}(\alpha)=\lambda_{0}(3\pi/4),\;\;\lim\limits_{\alpha\rightarrow 3\pi/4^+}\lambda_{1}(\alpha)=+\infty.\end{array}
\vspace{-0.2cm}$$

Let  $\mathcal{O}$ be a set in the space of the SLPs (1.1)--(1.2). Through the above observation, in the study of the continuity of the $n$-th eigenvalue function restricted in $\mathcal{O}$, it seems very important whether the number of eigenvalues of each  SLP in  $\mathcal{O}$ is equal or not.  In fact, we shall show that all
the eigenvalue functions restricted in $\mathcal{O}$ are continuous  if and only if the number of eigenvalues of  each  SLP
in $\mathcal{O}$  is equal (see Theorem 2.1).

Our study in the present paper was inspired by the remarkable work in [10] for the continuous case. Note that the
number of the eigenvalues  is finite for the discrete problem (1.1)--(1.2) and infinite for the continuous problem
(1.5)--(1.6). This difference will result in some differences between properties of their eigenvalue functions.
Consequently, our method used in the present paper is quite different from that used in the continuous case.
We shall list six aspects on these differences as follows. {\bf Firstly}, it was shown in [10] that if the minimal eigenvalue function $\lambda_0$  is bounded from below
 in a set of the space of the SLPs (1.5)--(1.6), then the $n$-th eigenvalue function restricted in  the set
 is continuous for each $n\geq0$. However, a similar conclusion is not true in the discrete case (see Example 2.1). Instead,
 as we have remarked in the above,  the continuity of the eigenvalue functions in a set can be completely determined by the number of eigenvalues of each SLP in the set (see Theorem 2.1).
 {\bf Secondly}, unlike that in the  continuous case, the $n$-th eigenvalue function
is not continuously dependent on the SLE (1.1) in general in the  discrete case.
{\bf Thirdly}, the discontinuity set of the $n$-th eigenvalue function in the space of BCs in the continuous case is different from that in the discrete case (see (1.8) and (4.12)).  They may be identified in a certain sense by letting $f_0\rightarrow+\infty$ in the discrete case.
{\bf Fourthly}, the authors in [10] employed  the Pr\"{u}fer transformation of (1.5) and some inequalities among eigenvalues of  SLPs (1.5)--(1.6) given in [7] in their discussions.
Though a discrete Pr\"{u}fer transformation was established in [2] and several inequalities among eigenvalues of discrete SLPs were obtained in [19, 20], we have found that it is quite difficult for us to employ a similar method to study this discrete problem.
Instead, we shall directly study several properties of the eigenvalue functions and make use of some spectral results of second-order difference equations given in [17, 22].
In order to  study  asymptotic behaviors of the $n$-th eigenvalue function $\lambda_n$   near a discontinuity point in the space of the SLPs (1.1)--(1.2), we shall first study  asymptotic  behaviors of  $\lambda_n$  in a certain direction near the discontinuity point. This way is shown to be convenient for our study in the discrete case.
{\bf Fifthly},  continuous dependence of the eigenvalue functions on the BCs is investigated via the local coordinate systems in the space of  BCs (1.2) directly, instead of being divided into  the separated and coupled cases.
{\bf Finally},  asymptotic behaviors of the $n$-th eigenvalue function near a discontinuity point in the space of BCs (1.2) are more complicated than those in the continuous case.

The rest of this paper is organized as follows. In Section 2, some notations and lemmas are introduced. A necessary and sufficient condition
 for all the  eigenvalue functions to be continuous and several
properties of the  eigenvalue functions in a set of the space of SLPs are given.   In Section 3,
  continuous  dependence of the $n$-th eigenvalue function on the SLE  is completely characterized for a fixed  BC.
In Section 4, continuous dependence of the $n$-th eigenvalue function on the BC is completely characterized  for a fixed SLE.
 In Section 5,  continuous dependence of the $n$-th eigenvalue function on the  SLP is studied.
 Then its continuity and discontinuity in the space of  the SLPs are completely characterized.\medskip

\noindent{\bf Remark 1.1.} We shall apply some results obtained in the present paper to study inequalities among the eigenvalues of general self-adjoint SLPs in our forthcoming paper.

\bigskip

\noindent{\bf 2. Preliminaries }\medskip

In this section, some notations and lemmas  are introduced. This section is divided into three parts.
In Section 2.1, the description of the space of the  SLPs is introduced.
Section 2.2 collects some basic properties of eigenvalues of the SLPs.
In Section 2.3,  a necessary and sufficient condition for all the eigenvalue functions to be continuous  and several
properties of the  eigenvalue functions in a set of the space of the SLPs are given. They are useful in the sequent sections.\medskip

\noindent{\bf 2.1. Space of self-adjoint discrete SLPs  }\medskip

In order to discuss continuous dependence of the $n$-th eigenvalue on the SLP, we need to know how to measure the closeness of two
SLEs and  of two BCs.

Let the SLE (1.1) be abbreviated as $(1/f,q,w)$.  Then the space of the SLEs can be written as \vspace{-0.1cm}$$\begin{array}{rrll}\Omega_N^{\mathbb{R},+} :=\{(1/f,q,w)\in\mathbb{R}^{3N+1}: {\rm (1.3)\; holds}\},\end{array} \vspace{-0.2cm}$$
and is equipped with the topology deduced from the real space $ \mathbb{R}^{3N+1}$. Note that $\Omega_N^{\mathbb{R},+}$ has $2^{N+1}$ connected components.
Bold faced lower case Greek letters, such as $ \pmb\omega$, are used to denote elements of $\Omega_N^{\mathbb{R},+}$.
For convenience, the maximum norm
in $ \mathbb{R}^{3N+1}$ will be used:
\vspace{-0.2cm}
$$\|(1/f,q,w)\|=\max\left\{|1/f_0|,\max_{1\leq n\leq N}\left\{|1/f_n|,|q_n|,|w_n|\right\}\right\}.\vspace{-0.2cm}$$

The quotient space
\vspace{-0.2cm}
$$\mathcal A^{\mathbb{C}} :=\raise2pt\hbox{${\rm M}^*_{2,4}(\mathbb C)$}/\lower3pt\hbox{${\rm
GL}(2,\mathbb C)$},
\vspace{-0.2cm}$$
equipped with the quotient topology, is taken as the space of general BCs; that is, each BC is an equivalence class of coefficient matrices of
system (1.2), where
\vspace{-0.1cm}$$\begin{array}{l}M_{2,4}^{*}(\mathbb{C}):=\{2\times4 \;{\rm complex\; matrix\; (A,B):}\;{\rm rank}(A,B)=2\},\\[1.0ex]
GL(2,\mathbb{C}):=\{2\times2\; {\rm comlplex \;matrix}\; T:{\rm det}\; T \neq0\}.\end{array}\vspace{-0.1cm} $$
The BC represented by (1.2) is denoted by $[A\,|\,B]$. Bold faced capital Latin letters, such as $\mathbf{A}$, are also used for BCs.
The space of self-adjoint BCs is denoted by $\mathcal{B}^\mathbb{C}$. The following result gives the topology and geometric structure of $\mathcal{B}^\mathbb{C}$.\medskip

\noindent{\bf Lemma 2.1} {\rm [22, Theorem 2.2]}. {\it The space $\mathcal{B}^\mathbb{C}$ equals the union of the following relative open sets:\vspace{-0.2cm}
$$\begin{array} {cccc} \mathcal{O}_{1,3}^{\mathbb{C}}=\left \{\left [\begin{array} {cccc}1&a_{12}&0&\bar{z}\\
0&z&-1&b_{22}\end{array}  \right ]:\; a_{12},b_{22}\in\mathbb{R},z\in \mathbb{C}\right \},\vspace{3mm}\\
 \mathcal{O}_{1,4}^{\mathbb{C}}=\left \{\left [\begin{array} {cccc}1&a_{12}&\bar{z}&0\\
0&z&b_{21}&1\end{array}  \right ]:\; a_{12},b_{21}\in\mathbb{R},z\in \mathbb{C}\right \},\vspace{3mm}\\
\mathcal{O}_{2,3}^{\mathbb{C}}=\left \{\left [\begin{array} {cccc}a_{11}&-1&0&\bar{z}\\
z&0&-1&b_{22}\end{array}  \right ]:\; a_{11},b_{22}\in\mathbb{R},z\in \mathbb{C}\right \},\vspace{3mm}\\
 \mathcal{O}_{2,4}^{\mathbb{C}}=\left \{\left [\begin{array} {cccc}a_{11}&-1&\bar{z}&0\\
z&0&b_{21}&1\end{array}  \right ]:\; a_{11},b_{21}\in\mathbb{R},z\in \mathbb{C}\right \}.\end{array}\vspace{-0.2cm}
                                                                                                                                                      \eqno(2.1)
$$
Moreover, $\mathcal{B}^\mathbb{C}$ is a connected and compact real-analytic manifold of dimension {\rm 4}.}\medskip

Lemma 2.1 says that $\mathcal{O}_{1,3}^{\mathbb{C}}$, $\mathcal{O}_{1,4}^{\mathbb{C}}$, $\mathcal{O}_{2,3}^{\mathbb{C}}$, and $\mathcal{O}_{2,4}^{\mathbb{C}}$ together form an atlas of local coordinate systems in $\mathcal{B}^\mathbb{C}$.
     If $\mathbf A\in\mathcal{O}_{1,3}^{\mathbb{C}}$, then the corresponding coordinate chart is given by
      \vspace{-0.2cm}$$\hat \varphi:\left [\begin{array} {cccc}1&a_{12}&0&\bar{z}\\
0&z&-1&b_{22}\end{array}  \right ]\to(a_{12},a,b,b_{22}), $$
      where $z=a+ib$, $a$, $b\in \mathbb R$.
      Others are given similarly. Thus, each of the four open sets in (2.1) can be identified with $\mathbb{R}^4$.

          The space $\Omega_N^{\mathbb{R},+}\times\mathcal{B}^\mathbb{C}$ of the SLPs is a real-analytic manifold of dimension
$3N+5$ and has $2^{N+1}$ connected components.

      The following result gives the canonical forms of separated and coupled self-adjoint BCs, respectively.\medskip

\noindent{{\bf Lemma 2.2} {\rm [21, Theorem 10.4.3]}. {\it Each separated self-adjoint BC can be written as \vspace{-0.1cm}
$$\mathbf{S}_{\alpha,\beta}:=\left [\begin{array} {cccc}\cos\alpha&-\sin\alpha&0&0\\
0&0&\cos\beta&-\sin\beta\end{array}  \right ],                                                                                                        \eqno(2.2)
\vspace{-0.2cm}$$
where \vspace{-0.1cm}$$\alpha\in[0,\pi),\beta\in (0,\pi];\vspace{-0.1cm}$$
and each coupled self-adjoint BC can be written as \vspace{-0.2cm}
$$[e^{i\gamma}K\,|\,-I],\vspace{-0.2cm}$$
where\vspace{-0.2cm}
$$\gamma\in[0,\pi),\;\;K\in SL(2,\mathbb{R}):=\{2\times2 \;real\; matrix\; M:{\rm det}M=1\}.\vspace{-0.2cm}$$}\vspace{-0.2cm}

The spaces of separated and coupled self-adjoint BCs are denoted by $\mathcal{B}_S$ and  $\mathcal{B}_C$, respectively. So $\mathcal{B}^{\mathbb{C}}=\mathcal{B}_S\cup\mathcal{B}_C$.

\medskip

\noindent{\bf 2.2. Basic properties of eigenvalues  }\medskip

In this subsection, some basic properties of eigenvalues of the SLPs are introduced.\medskip

For each $\lambda\in\mathbb{C}$, let $\phi(\lambda)$ and $\psi(\lambda)$ be the solutions of (1.1)
satisfying the following initial conditions: \vspace{-0.2cm}
$$
\phi_0(\lambda)=1, f_0\Delta\phi_0(\lambda)=0; \;\;
\psi_0(\lambda)=0, \ f_0\Delta\psi_0(\lambda)=1,
\vspace{-0.2cm}$$
separately. It follows from [22] that the leading terms of $\phi_N (\lambda)$, $\psi_N (\lambda)$, $f_N\Delta\phi_N (\lambda)$,
and $f_N\Delta\psi_N (\lambda)$ as polynomials of $\lambda$ are \vspace{-0.2cm} $$\begin{array}{cccc}(-1)^{N-1}\left(\prod\limits_{i=1}^{N-1}({w_i}/{f_i})\right)
\lambda^{N-1},&(-1)^{N-1}\left(({1}/{f_0})\prod\limits_{i=1}^{N-1}({w_i}/{f_i})\right)
\lambda^{N-1},
\\[2.0ex](-1)^{N}\left(w_N\prod\limits_{i=1}^{N-1}({w_i}/{f_i})\right)
\lambda^{N},&(-1)^{N}\left(({w_N}/{f_0})\prod\limits_{i=1}^{N-1}({w_i}/{f_i})\right)
\lambda^{N},\end{array}                                                                                                                               \eqno(2.3)
\vspace{-0.2cm} $$
respectively.
\medskip

\noindent{\bf Lemma 2.3 } {\rm [22, Lemmas 3.2 and 3.3]}. {\it A number $\lambda \in
\mathbb{R}$ is an eigenvalue of the  SLP
{\rm(1.1)}-{\rm(1.2)}
 if and only if $\lambda$ is a zero of the characteristic polynomial\vspace{-0.2cm}
$$\Gamma(\lambda) =\det A +\det B +G(\lambda),\vspace{-0.2cm}$$
where}\vspace{-0.2cm}
$$G(\lambda):=c_{11}\phi_N(\lambda)+c_{12}\psi_N(\lambda)+c_{21}f_N\Delta\phi_N(\lambda)
+c_{22}f_N\Delta\psi_N(\lambda),\vspace{-0.2cm}$$
\vspace{-0.2cm}$$ C: =\left(\begin{array}{cc} b_{11} & b_{21} \\ b_{12} &b_{22}\end{array}\right)
 \left(\begin{array}{cc} a_{22} & -a_{21} \\ -a_{12} & a_{11} \end{array}\right).
\vspace{-0.2cm}$$
 \medskip

When we count the
 eigenvalues of an SLP in a domain in $\mathbb{R}$, their multiplicities will be
taken into account.

Let $(\pmb\omega,\mathbf{A})\in\Omega_N^{\mathbb{R},+} \times \mathcal{B}^{\mathbb{C}}$. Set \vspace{-0.2cm}
$$r=r(\pmb\omega,\mathbf{A}):={\rm rank} \left(\begin{array}{cc}-a_{11}+f_0a_{12}&b_{12}\\-a_{21}+f_0a_{22}&b_{22}\end{array}\right).\eqno(2.4)\vspace{-0.2cm}
$$
 Obviously, $0\leq r \leq 2$.\medskip

\noindent{\bf Lemma 2.4} [22, Lemma 3.4]. {\it
The number of eigenvalues of $(\pmb\omega,\mathbf{A})$  is equal to $N-2+r$, where $r$ is defined by {\rm(2.4)}. }\medskip

The above results can be  deduced from [17, Theorem 4.1].

 By (2.3) and Lemma 2.3, the coefficient of $\lambda^N$ in the polynomial $\Gamma(\lambda)$ is \vspace{-0.2cm}
 $$\theta(\pmb\omega,\mathbf{A}):=(-1)^{N}\left({w_N}\prod\limits_{i=1}^{N-1}({w_i}/{f_i})\right)\left[(a_{11}b_{22}-a_{21}b_{12})
/f_0+a_{22}b_{12}-a_{12}b_{22}\right].                                                                                                           \eqno(2.5)
 \vspace{-0.4cm}$$
 Thus, by Lemma 2.4,  there are exactly the following  three  cases:
 \begin{itemize}\vspace{-0.2cm}
\item[{\rm (i)}]  $(\pmb\omega,\mathbf {A})$ has exactly $N$ eigenvalues in the case that  $\theta(\pmb\omega,\mathbf {A})\neq0$;\vspace{-0.2cm}
\item[{\rm (ii)}] $(\pmb\omega,\mathbf {A})$ has exactly $N-1$ eigenvalues in the case that  $\theta(\pmb\omega,\mathbf {A})=0$ and $\mathbf {A}\neq\mathbf {A}(1/f_0)$, where \vspace{-0.2cm}$$\mathbf {A}(1/f_0)=\left [\begin{array} {cccc}1&1/f_0&0&0\\
0&0&1&0\end{array}  \right ];\vspace{-0.4cm}$$ \vspace{-0.5cm}
\item[{\rm (iii)}] $(\pmb\omega,\mathbf {A})$ has exactly $N-2$ eigenvalues in the case that $\theta(\pmb\omega,\mathbf {A})=0$ and $\mathbf {A}=\mathbf {A}(1/f_0)$.
 \vspace{-0.5cm}\end{itemize}\vspace{-0.4cm}

By $\lambda_n(\pmb\omega,\mathbf {A})$ denotes the $n$-th eigenvalue of $(\pmb\omega,\mathbf {A})$. When the SLE is fixed, $\lambda_n(\mathbf{A})$ is also used  for $\mathbf{A}\in \mathcal{B}^{\mathbb{C}}$; when the BC is fixed,
$\lambda_n(\pmb\omega)$ is also used for $\pmb\omega\in \Omega_N^{\mathbb{R},+}$, and etc.

The following result is a generation of  [22, Corollary 3.3].\medskip

\noindent{\bf Lemma 2.5.} {\it Assume that $\mathcal{O}$ is a set of $\Omega_N^{\mathbb{R},+} \times \mathcal{B}^{\mathbb{C}}$, and
 $(\pmb\omega_0,\mathbf {A}_0) \in\mathcal{O}$. Let $r_1$ and $r_2$ be two real numbers with  $r_1<r_2$ such that neither
of them is an eigenvalue of $(\pmb\omega_0,\mathbf {A}_0)$, and $n\geq 0$ be the number of eigenvalues of $(\pmb\omega_0,\mathbf {A}_0)$
in the interval $[r_1,r_2]$. Then there exists a neighborhood $\mathcal{U}$ of $(\pmb\omega_0,\mathbf {A}_0)$ in
$\mathcal{O}$ such that each $(\pmb\omega,\mathbf {A})\in\mathcal{U}$  has
 exactly $n$ eigenvalues in $[r_1,r_2]$, which all lie in $(r_1,r_2)$.}\medskip

\noindent{\bf Proof.} Since the proof is similar to that of  [22, Theorem 3.4], we omit its details.\medskip

\noindent{\bf Lemma 2.6} {\rm [22, Theorem 3.5 and Remark 3.2]}.
 {\it Assume that $\mathcal{O}$ is a connected set of $\Omega_N^{\mathbb{R},+} \times \mathcal{B}^{\mathbb{C}}$ and
 $(\pmb\omega_0,\mathbf {A}_0) \in\mathcal{O}$. Let $\lambda_*$ be an eigenvalue of $(\pmb\omega_0,\mathbf A_0)$ with multiplicity $m$.
Fix a small $\epsilon>0$ such that $\lambda_*$ is the only
eigenvalue of $(\pmb\omega_0,\mathbf A_0)$ in the interval
$[\lambda_*-\epsilon, \lambda_*+\epsilon]$. Then there is a connected neighborhood $\mathcal{F}\subset\mathcal{O}$ and continuous functions $\Lambda_i:\mathcal F \to \mathbb R$, $1\leq i\leq m$, such that}\vspace{-0.2cm}
\begin{itemize}
\item[{\rm (i)}] {\it$\Lambda_i(\pmb\omega_0,\mathbf A_0)
=\lambda_*, 1\leq i \leq m$;}\vspace{-0.2cm}
\item[{\rm (ii)}] {\it$\lambda_*-\epsilon <\Lambda_1(\pmb\omega,\mathbf A) \leq\cdots\leq
\Lambda_m(\pmb\omega,\mathbf A) <\lambda_*+\epsilon$ for each
$(\pmb\omega,\mathbf A) \in \mathcal F$;}\vspace{-0.2cm}
\item[{\rm (iii)}] {\it
$\Lambda_i(\pmb\omega,\mathbf A), 1\leq i \leq m$, are eigenvalues of
$(\pmb\omega,\mathbf A)$ for each $(\pmb\omega,\mathbf A) \in \mathcal F$.}
\end{itemize}
\noindent{\bf Remark 2.1.} In Lemma 2.6, $\Lambda_i, 1\leq i\leq m,$  are  called the continuous eigenvalue
branches through
$\lambda_*$.\medskip

\noindent{\bf 2.3. Properties of the eigenvalue functions}\medskip

In this subsection, a necessary and sufficient condition  for all the eigenvalue functions  to be continuous
and several other properties  of the eigenvalue functions in a set of the space of the SLPs, which  are useful  in the study of asymptotic behaviors of the $n$-th eigenvalue function near a discontinuity point, are obtained.  \medskip

 \noindent{\bf Theorem 2.1.} {\it Let $\mathcal{O}$ be a  set of $\Omega_N^{\mathbb{R,+}} \times \mathcal{B}^{\mathbb{C}}$.
Then the number of eigenvalues of each $(\pmb\omega,\mathbf {A})\in\mathcal{O}$ is equal if and only if
each eigenvalue function $\lambda_n(\pmb\omega,\mathbf {A})$ restricted in $\mathcal{O}$  is
 continuous. Furthermore, if  $\mathcal{O}$ is a connected set of $\Omega_N^{\mathbb{R,+}} \times \mathcal{B}^{\mathbb{C}}$, then each eigenvalue function $\lambda_n(\pmb\omega,\mathbf {A})$ is locally
a continuous eigenvalue branch in $\mathcal{O}$.}\medskip

\noindent{\bf Proof.} Suppose that each $(\pmb\omega,\mathbf {A})\in\mathcal{O}$ has exactly $k$ eigenvalues: $\lambda_0(\pmb\omega,\mathbf {A})\leq\lambda_1(\pmb\omega,\mathbf {A})\leq\cdots\leq\lambda_{k-1}(\pmb\omega,\mathbf {A})$. Then the $k$ eigenvalue functions $\lambda_i :\mathcal{O}\rightarrow\mathbb{R}, 0\leq i\leq k-1,$
are well-defined. Now, we show that $\lambda_i$ restricted in $\mathcal{O}$ is  continuous at a fixed $(\pmb\omega_0,\mathbf {A}_0) \in\mathcal{O}$.
We first consider the case that  $\lambda_i(\pmb\omega_0,\mathbf {A}_0)$ is a simple eigenvalue.
Fix $0< i< k-1$. Let $r_j, 1\leq j \leq4$,  be four real numbers such that
  $r_1<\lambda_{0}(\pmb\omega_0,\mathbf {A}_0)$, $\lambda_{i-1}
(\pmb\omega_0,\mathbf {A}_0)<r_2<\lambda_i(\pmb\omega_0,\mathbf {A}_0)<r_3<\lambda_{i+1}(\pmb\omega_0,
\mathbf {A}_0)$ and $r_4>\lambda_{k-1}(\pmb\omega_0,\mathbf {A}_0)$. By Lemma 2.5 there exists a neighborhood
  $\mathcal{U}$  of $(\pmb\omega_0,\mathbf {A}_0)$  in $\mathcal{O}$ such that each $(\pmb\omega,\mathbf {A})\in\mathcal{U}$ has
 exactly $i$ eigenvalues in $(r_1,r_2)$, exactly one eigenvalue in $(r_2,r_3)$, and
   exactly $k-i-1$ eigenvalues in $(r_3,r_4)$.
   Since each $(\pmb\omega,\mathbf {A})\in\mathcal{U}$ has exactly $k$ eigenvalues,
the eigenvalue of $(\pmb\omega,\mathbf {A})$ in $(r_2,r_3)$ is exactly $\lambda_i(\pmb\omega,\mathbf {A})$.
 Hence,  $\lambda_i$ restricted in $\mathcal{O}$ is continuous at $(\pmb\omega_0,\mathbf {A}_0)$.
 With a similar method used above, one can easily verify that $\lambda_{0}$ and $\lambda_{k-1}$ restricted in $\mathcal{O}$ are continuous at $(\pmb\omega_0,\mathbf {A}_0)$.

Suppose that $\mathcal{O}$ is connected. Then the above $\mathcal{U}$ can be chosen  to be also connected. By Lemma 2.6 and Remark 2.1,
 $\lambda_i$ restricted in $\mathcal{U}$ is exactly a continuous eigenvalue branch  through $\lambda_i(\pmb\omega_0,\mathbf {A}_0)$.
Thus, $\lambda_i$  is locally a continuous eigenvalue branch in $\mathcal{O}$ for each $0\leq i\leq k-1$.

 In the case that  the multiplicity of $\lambda_i(\pmb\omega_0,\mathbf {A}_0)$ is equal to 2,
one can show that the results still hold with a similar argument.

Conversely, suppose that there exists an $(\pmb\omega_1,\mathbf {A}_1)\in\mathcal{O}$ such that   the number of  eigenvalues of $(\pmb\omega_1,\mathbf {A}_1)$    is not equal to that of another. Then there exists an
eigenvalue function $\lambda_{i_0}$ that   can not be well-defined at least at one point in $\mathcal{O}$. Thus $\lambda_{i_0}$ restricted in $\mathcal{O}$ is not
continuous.  This completes the proof.\medskip

\noindent{\bf Remark 2.2.} By Theorem 2.1, if the number of eigenvalues of each $(\pmb\omega,\mathbf {A})\in\mathcal{O}$ is not equal, then there exists
at least one eigenvalue function that is not continuous in $\mathcal{O}$. In this case, are all the eigenvalue functions not continuous in $\mathcal{O}$? The following example gives a negative answer.\medskip

 \noindent{\bf Example 2.1.} \rm Let $s\in[0,2]$. Consider the 1-parameter family of
 SLPs (1.1)--(1.2) with
\vspace{-0.1cm}$$\begin{array}{cccc}
f_0 =\begin{cases} {1\over 2-s}  & \text{ if } s\in [0,1), \\
                     {1\over s} & \text{ if } s\in[1,2], \end{cases} \qquad
f_1 =\begin{cases} 1 & \text{ if } s\in [0,1), \\
                     {1\over s} & \text{ if } s\in[1,2],\end{cases}\\[5.0ex]
\ f_2=1, \; q_1 =q_2 =0, \; w_1 =w_2 =1,\; N=2,
\vspace{-0.2cm}\end{array}$$
and
\vspace{-0.2cm}$$
\mathbf A =\left[\begin{array}{cccc} 1 & 1 & 0 & 0 \\ 0 & 0 & -1 & 1 \end{array}\right].
                                                                                                                                             \eqno(2.6)
\vspace{-0.2cm}$$
Then, by Lemma 2.3, direct calculations deduce that the
characteristic function is
$$\vspace{-0.2cm}
\Gamma(\lambda) =\begin{cases} (1-s)\lambda^2+(s-2)\lambda+s-1 & \text{ if } s\in [0,1), \\
                     (s^2-s)\lambda^2+(s^2-4s+2)\lambda+2-2s & \text{ if } s\in[1,2]. \end{cases}
\vspace{-0.05cm}$$
Thus, the SLP with $s=1$ has exactly  one eigenvalue and with each $s\in[0,1)\cup(1,2]$ has exactly  two eigenvalues, which are given as
\vspace{-0.1cm}$$\begin{array}{cccc}
\lambda_0(s) =\begin{cases} {2-s-\sqrt{5s^2-12s+8} \over 2(1-s)}  & \text{ if } s\in [0,1), \\
                            0 & \text{ if } s=1,\\
                     {-s^2+4s-2 -\sqrt{(s^2-4s+2)^2-4(s^2-s)(2-2s)} \over
  2(s^2-s) } & \text{ if } s\in (1,2],\end{cases} \qquad\\[5.0ex]
\lambda_1(s) =\begin{cases} {2-s+\sqrt{5s^2-12s+8} \over 2(1-s)}& \text{ if } s\in [0,1), \\
                     {-s^2+4s-2 +\sqrt{(s^2-4s+2)^2-4(s^2-s)(2-2s)} \over
  2(s^2-s) }  & \text{ if }  s\in (1,2]. \end{cases}
\end{array}\vspace{-0.2cm}$$
\vspace{-0.5cm}
\begin{center}

\includegraphics[width=87mm]{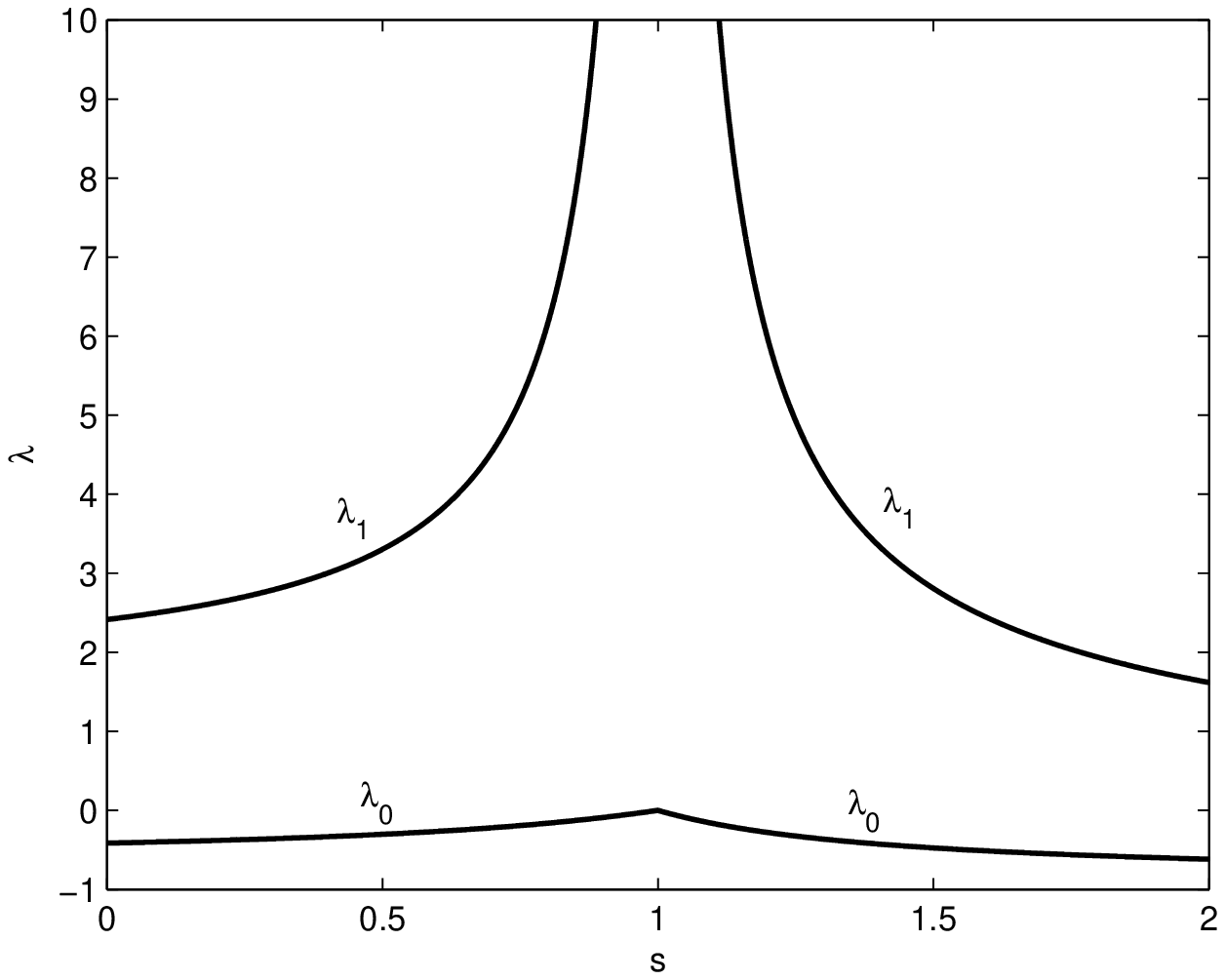}

{\small {\bf Figure 2.1.}}
\end{center}
\vspace{-0.2cm}
See Figure 2.1. Then the number of eigenvalues of the SLP  with $s=1$ is not equal to that of another one. It is evident that the eigenvalue function $\lambda_0$ is continuous in $s\in[0,2]$.  Note that   $\lambda_0$
 is bounded from below for $s\in[0,2]$. However, the eigenvalue function $\lambda_1$ is not well-defined at $s=1$, and thus   $\lambda_1$ is not continuous at $s=1$.\medskip

Now we give several other properties of the eigenvalue functions, which are useful in the discussions
 in the sequent sections.
\medskip

\noindent{\bf Theorem 2.2.} {\it Assume that $\mathcal{O}\subset\Omega_N^{\mathbb{R,+}} \times \mathcal{B}^{\mathbb{C}}$ satisfies that each $(\pmb\omega,\mathbf {A})$ in $\mathcal{O}$ has exactly $k$ eigenvalues, where $k\geq2$. Let $(\pmb\omega_0,\mathbf {A}_0) \in\bar{\mathcal{O}}\backslash\mathcal{O}$ have exactly $m$ eigenvalues for some $0< m\leq k-1$.
 \begin{itemize}
\item[{\rm (i)}]\vspace{-0.2cm}
If the first $k-m$ eigenvalue functions satisfy
\vspace{-0.2cm}$$\lim\limits_{\mathcal{O}\ni(\pmb\omega,\mathbf {A})\rightarrow(\pmb\omega_0,\mathbf {A}_0)}\lambda_n(\pmb\omega,\mathbf {A})=-\infty,\;\;
0\leq n\leq k-m-1,
                                                                                                                       \eqno(2.7)\vspace{-0.2cm}$$
then the last $m$ eigenvalue functions satisfy
\vspace{-0.2cm}$$\lim\limits_{\mathcal{O}\ni(\pmb\omega,\mathbf {A})\rightarrow(\pmb\omega_0,\mathbf {A}_0)}\lambda_n(\pmb\omega,\mathbf {A})=
\lambda_{n-k+m}(\pmb\omega_0,\mathbf {A}_0),\;\;k-m\leq n \leq k-1.                                                                           \eqno(2.8)
\vspace{-0.20cm}$$
\vspace{-1cm}
\item[{\rm (ii)}] If the last $k-m$ eigenvalue functions satisfy
\vspace{-0.2cm}$$\lim\limits_{\mathcal{O}\ni(\pmb\omega,\mathbf {A})\rightarrow(\pmb\omega_0,\mathbf {A}_0)}\lambda_n(\pmb\omega,\mathbf {A})=+\infty
,\;\;m\leq n\leq k-1,\vspace{-0.2cm}$$
then the first $m$ eigenvalue functions satisfy
\vspace{-0.2cm}$$\lim\limits_{\mathcal{O}\ni(\pmb\omega,\mathbf {A})\rightarrow(\pmb\omega_0,\mathbf {A}_0)}
\lambda_n(\pmb\omega,\mathbf {A})=\lambda_{n}(\pmb\omega_0,\mathbf {A}_0),\;\;0\leq n \leq m-1.\vspace{-0.2cm}$$
\end{itemize}}

\noindent{\bf Proof.} First, we show that (i) holds. Let $r_1<r_2$ be any two real numbers such that $r_1<\lambda_0(\pmb\omega_0,\mathbf {A}_0)$ and $r_2>\lambda_{m-1}(\pmb\omega_0,\mathbf {A}_0)$.
It follows from (2.7) that there exists a neighborhood $\mathcal{U}$  of $(\pmb\omega_0,\mathbf {A}_0)$ in $ \mathcal{O}\cup\{(\pmb\omega_0,\mathbf {A}_0)\}$
such that
\vspace{-0.2cm}$$\lambda_n(\pmb\omega,\mathbf {A})<r_1,\;\forall\;(\pmb\omega,\mathbf {A})\in\mathcal{U}\backslash\{(\pmb\omega_0,\mathbf {A}_0)\},\; 0\leq n\leq k-m-1.                                                      \eqno(2.9)
\vspace{-0.5cm}$$
Since $(\pmb\omega_0,\mathbf {A}_0)$ has exactly $m$ eigenvalues, by Lemma 2.5 there exists a neighborhood $\mathcal{U}_1\subset\mathcal{U}$
of $(\pmb\omega_0,\mathbf {A}_0)$ such that each $(\pmb\omega,\mathbf {A})\in\mathcal{U}_1$ has exactly $m$ eigenvalues  in $[r_1,r_2]$, which are all in
$(r_1,r_2)$. Hence,
for each $(\pmb\omega,\mathbf {A})\in\mathcal{U}_1\backslash\{(\pmb\omega_0,\mathbf {A}_0)\}$, by (2.9) and noting that $(\pmb\omega,\mathbf {A})$
 has exactly $k$ eigenvalues, one has that
 \vspace{-0.2cm}$$\lambda_j(\pmb\omega,\mathbf {A})\in(r_1,r_2),\;\; k-m\leq j\leq k-1.                                                     \eqno(2.10)
\vspace{-0.2cm}$$

Now, suppose that $\lambda_i(\pmb\omega_0,\mathbf {A}_0)$ is a simple eigenvalue for some $0\leq i\leq m-1$.
Let $(r_3,r_4)\subset(r_1,r_2)$ be an interval such that $\lambda_i(\pmb\omega_0,\mathbf{A}_0)$ is the only eigenvalue of $(\pmb\omega_0,\mathbf{A}_0)$ in
$[r_3,r_4]$, and $\lambda_i(\pmb\omega_0,\mathbf{A}_0)\in(r_3,r_4)$.
Again, by Lemma 2.5 there exists a neighborhood
  $\mathcal{U}_2\subset\mathcal{U}_1$ of $(\pmb\omega_0,\mathbf {A}_0)$ such that each $(\pmb\omega,\mathbf {A})\in\mathcal{U}_2$
has exactly $i$ eigenvalues in $(r_1,r_3)$, exactly one eigenvalue in $(r_3,r_4)$, and
   exactly $m-i-1$ eigenvalues in $(r_4,r_2)$.
 Thus, it follows from (2.9) and (2.10)  that $\lambda_{k-m+i}(\pmb\omega,\mathbf {A})\in(r_3,r_4)$ for each $(\pmb\omega,\mathbf {A})\in\mathcal{U}_2\backslash\{(\pmb\omega_0,\mathbf {A}_0)\}$.  Consequently,
  \vspace{-0.2cm}$$\lim\limits_{\mathcal{O}\ni(\pmb\omega,\mathbf {A})\rightarrow(\pmb\omega_0,\mathbf {A}_0)}\lambda_{k-m+i}(\pmb\omega,\mathbf {A})=
\lambda_{i}(\pmb\omega_0,\mathbf {A}_0). \vspace{-0.2cm}$$
 Therefore, (2.8) holds in this case.

In the other case that the multiplicity of  $\lambda_i(\pmb\omega_0,\mathbf {A}_0)$ is equal to 2
 for some $0\leq i\leq m-1$, one can show that the results hold with a similar argument.

Since the proof of (ii) is similar to that of (i), we omit its details.  The  proof is complete.\medskip

\noindent{\bf Theorem 2.3.} {\it Let $\mathcal{O}$ be a connected set of $\Omega_N^{\mathbb R,+}\times\mathcal{B}^{\mathbb{C}}$ and
$r_1$ and $r_2$ be two real numbers with $r_1<r_2$.
Assume that  each problem in $\mathcal{O} $ has exactly $k$ eigenvalues,
and exactly  $m\geq0$ eigenvalues in $[r_1,r_2]$ that are
 in $(r_1,r_2)$ with $m<k$. Then its other $k-m$ eigenvalues out of $[r_1,r_2]$, denoted by $\hat\lambda_1(\pmb\omega,\mathbf{A})\leq
 \cdots\leq\hat\lambda_{k-m}(\pmb\omega,\mathbf{A})$, have the following properties:
\begin{itemize}\vspace{-0.2cm}
\item[{\rm (i)}] For each $1\leq i \leq k-m$,
 \vspace{-0.2cm}$$ {  either}\; E_i:=\{\hat\lambda_i(\pmb\omega,\mathbf{A}):(\pmb\omega,\mathbf{A})\in\mathcal{O}\}\subset(-\infty,r_1)\; { or} \;E_i\subset(r_2,+\infty);
                                                                                                                                                     \eqno(2.11)
  \vspace{-0.2cm}$$
there exists $0\leq i_0\leq k-1$ such that $\hat\lambda_i(\pmb\omega,\mathbf{A})=\lambda_{i_0}(\pmb\omega,\mathbf{A})$ for all $(\pmb\omega,\mathbf{A})\in\mathcal{O}$, and consequently, $\hat\lambda_i$  is continuous in $\mathcal{O}$.\vspace{-0.2cm}
\item[{\rm (ii)}] Further, assume that $(\pmb\omega_0,\mathbf{A}_0)\in \bar{\mathcal{O}}\backslash\mathcal{O}$ has exactly $m$ eigenvalues
$\lambda_j(\pmb\omega_0,\mathbf{A}_0)$, $0\leq j\leq m-1$, and all of them are in $(r_1,r_2)$. If $E_i\subset(-\infty,r_1)$ for some $1\leq i\leq k-m$,  then
\vspace{-0.2cm}$$\lim\limits_{\mathcal{O}\ni(\pmb\omega,\mathbf{A})\rightarrow(\pmb\omega_0,\mathbf{A}_0)}\hat\lambda_i(\pmb\omega,\mathbf{A})=-\infty.
\vspace{-0.1cm}$$
If $E_i\subset(r_2,+\infty)$ for some $1\leq i\leq k-m$, then
\vspace{-0.2cm}$$\lim\limits_{\mathcal{O}\ni(\pmb\omega,\mathbf{A})\rightarrow(\pmb\omega_0,\mathbf{A}_0)}\hat\lambda_i(\pmb\omega,\mathbf{A})=+\infty.
\vspace{-0.2cm}$$
\end{itemize}}\vspace{-0.2cm}

\noindent{\bf Proof.} By the assumption that  each $(\pmb\omega,\mathbf{A})\in\mathcal{O}$ has exactly $k$ eigenvalues,
the $n$-th eigenvalue functions $\lambda_j$, $0\leq j \leq k-1$, restricted in $\mathcal{O}$  are continuous  by Theorem 2.1.

First, we show  that (i) holds.   We claim that (2.11) holds for $i=k-m$.
Otherwise, there exist $(\pmb\omega_1,\mathbf{A}_1)$, $(\pmb\omega_2,\mathbf{A}_2)\in\mathcal{O}$ such that
 $\hat\lambda_{k-m}(\pmb\omega_1,\mathbf{A}_1)\in(-\infty,r_1)$ and $\hat\lambda_{k-m}(\pmb\omega_2,\mathbf{A}_2)\in(r_2,+\infty)$.
 Then $\lambda_{k-1}(\pmb\omega_1,\mathbf{A}_1)\in(r_1,r_2)$ and $\lambda_{k-1}(\pmb\omega_2,\mathbf{A}_2)\in(r_2,+\infty)$.
 Since $r_2$ is not an eigenvalue of any problem in $\mathcal{O}$,
$\lambda_{k-1}$ is not continuous in $\mathcal{O}$.
This is a contradiction by the assumption that $\mathcal{O}$ is connected.

 If $E_{k-m}\subset(-\infty,r_1)$,
then $E_{i}\subset(-\infty,r_1)$ for all $1\leq i \leq k-m-1$,
  and consequently, $\hat \lambda_j(\pmb\omega,\mathbf{A})=\lambda_{j-1}(\pmb\omega,\mathbf{A})$
for all $(\pmb\omega,\mathbf{A})\in\mathcal{O}$ and $1\leq j\leq k-m$.

 If $E_{k-m}\subset(r_2,+\infty)$,
then using the same method as employed in the above paragraph, one can show that either
 $E_{k-m-1}\subset(-\infty,r_1)$ or $E_{k-m-1}\subset(r_2,+\infty)$.
If $E_{k-m-1}\subset(-\infty,r_1)$,
 then $E_{j}\subset(-\infty,r_1)$, $1\leq j \leq k-m-2$.
 Hence,  $\hat \lambda_{k-m}(\pmb\omega,\mathbf{A})=\lambda_{k-1}(\pmb\omega,\mathbf{A})$ and $\hat \lambda_j(\pmb\omega,\mathbf{A})=\lambda_{j-1}(\pmb\omega,\mathbf{A})$, $1\leq j\leq k-m-1$, for  all
 $(\pmb\omega,\mathbf{A})\in\mathcal{O}$.
 If $E_{k-m-1}\subset(r_2,+\infty)$,
then again using the same method as employed in the above discussion, one can show that
 either $E_{k-m-2} \subset(-\infty,r_1)$ or $E_{k-m-2} \subset(r_2,+\infty)$.
This procedure
can be finished in finite steps.

Further, since the $n$-th eigenvalue functions $\lambda_j,0\leq j \leq k-1$, restricted in $\mathcal{O}$  are continuous,
 $\hat \lambda_i$ is continuous in $\mathcal{O}$ for each $1\leq i \leq k-m$.  Then  (i) has been shown.

Now, we show that (ii) holds. We only consider the case that $E_{i}\subset(-\infty,r_1)$ for some $1\leq i\leq k-m$.
The other case can be  similarly discussed. Otherwise, there exists a positive number $M>|r_1|$ such that for any neighborhood $\mathcal{ U}$
 of $(\pmb\omega_0,\mathbf{A}_0)$  in $\mathcal{O}\cup\{(\pmb\omega_0,\mathbf{A}_0)\}$,
there exists $(\hat{\pmb\omega},\hat{\mathbf{A}})\in\mathcal{ U}\backslash\{(\pmb\omega_0,\mathbf{A}_0)\}$ satisfying
  $ -M\leq\hat\lambda_i(\hat{\pmb\omega},\hat{\mathbf{A}})<r_1$. Since $(\pmb\omega_0,\mathbf{A}_0)$ has exactly $m$ eigenvalues in $(-M, r_2)$,
 by Lemma 2.5 there exists a neighborhood $\mathcal{U}_1\subset\mathcal{O}\cup\{(\pmb\omega_0,\mathbf{A}_0)\}$ of $(\pmb\omega_0,\mathbf{A}_0)$
 such that each $(\pmb\omega,\mathbf{A})\in \mathcal{U}_1$ has exactly  $m$ eigenvalues in $(-M, r_2)$. However,
taking $\mathcal{U}\subset\mathcal{U}_1$, we get that $(\hat{\pmb\omega},\hat{\mathbf{A}})$ has
 at least $m+1$ eigenvalues in $(-M,r_2)$ by the assumption that each  $(\pmb\omega,\mathbf{A})\in\mathcal{O}$ has $m$ eigenvalues in $(r_1,r_2)$. This is a contradiction. The entire proof is complete.\medskip

In some cases, $(\pmb\omega,\mathbf{A})\in\Omega_N^{\mathbb R,+}\times\mathcal{B}^\mathbb{C}$ may be  continuously dependent on
some real parameters or variables. For example, $(\pmb\omega,\mathbf{A})$ is continuously dependent on $1/f_0$; $(\pmb\omega,
\mathbf{S}_{\alpha,\beta})$ is continuously dependent on $\alpha$ and $\beta$; $(\pmb\omega,\mathbf{A})$, where $\mathbf{A}
\in\mathcal{O}^{\mathbb{C}}_{1,3}$, is continuously dependent on variables $a_{12}$ and $b_{22}$, etc. Then, we  write
 $(\pmb\omega,\mathbf{A})_{\nu}$ instead of $(\pmb\omega,\mathbf{A})$ to indicate the  dependence of $(\pmb\omega,\mathbf{A})$
on a variable or parameter $\nu$ in some situations for convenience. Next, we shall discuss the dependence of $\lambda_n(\nu):=\lambda_n((\pmb\omega,\mathbf{A})_{\nu})$ on $\nu$.\medskip

\noindent{\bf Lemma 2.7.} {\it Let $(\pmb\omega,\mathbf{A})_{\nu}\in\Omega_N^{\mathbb R,+}\times\mathcal{B}^\mathbb{C}$ be continuously dependent on a real variable or parameter $\nu$ in $(\nu_0-\epsilon,\nu_0+\epsilon)$ for some $\epsilon>0$ and  $\mathcal{O}:=\{(\pmb\omega,\mathbf{A})_{\nu}:\nu\in(\nu_0-\epsilon,
\nu_0+\epsilon)\}$.  Assume that $(\pmb\omega,\mathbf{A})_{\nu_0}$ has exactly $m\geq1$ eigenvalues and for each $\nu\in(\nu_0-\epsilon,\nu_0+\epsilon)\backslash \{\nu_0\}$, $(\pmb\omega,\mathbf{A})_{\nu}$  has exactly $k$ eigenvalues with $k>m$.

\begin{itemize}\vspace{-0.2cm}
\item[{\rm (i)}] If the $n$-th eigenvalue functions $\lambda_n(\nu)$ are non-increasing in $(\nu_0-\epsilon_0,\nu_0)$ for all $0\leq n \leq k-1$, then they  satisfy the following asymptotic behaviors near $\nu_0$:
    \vspace{-0.2cm} $$\lim\limits_{\nu\rightarrow\nu_0^-}\lambda_n({\nu})=-\infty,\; 0\leq n \leq k-m-1,         \eqno(2.12) \vspace{-0.7cm} $$
     \vspace{-0.7cm} $$
  \lim\limits_{\nu\rightarrow\nu_0^-}\lambda_n({\nu})=\lambda_{n-k+m}({\nu_0}), \; k-m\leq n \leq k-1. \eqno(2.13)
  \vspace{-0.9cm} $$  \vspace{-0.7cm}
  \item[{\rm (ii)}] If the $n$-th eigenvalue functions $\lambda_n({\nu})$ are non-decreasing
 in $(\nu_0-\epsilon_0,\nu_0)$  for all $0\leq n \leq k-1$, then
they
 satisfy the following asymptotic behaviors near $\nu_0$:
    \vspace{-0.2cm}$$
    \lim\limits_{\nu\rightarrow\nu_0^-}\lambda_n({\nu})=\lambda_{n}({\nu_0}), \;\; 0\leq n \leq m-1,\;
     \lim\limits_{\nu\rightarrow\nu_0^-}\lambda_n({\nu})=+\infty,\; m\leq n \leq k-1.
    \vspace{-0.5cm}$$
   \item[{\rm (iii)}] If the $n$-th eigenvalue functions $\lambda_n(\nu)$ are non-increasing in $ (\nu_0,\nu_0+\epsilon_0)$ for all $0\leq n \leq k-1$, then they  satisfy the following asymptotic behaviors near $\nu_0$:
  \vspace{-0.2cm} $$
\lim\limits_{\nu\rightarrow\nu_0^+}\lambda_n({\nu})=\lambda_{n}({\nu_0}), \;\; 0\leq n \leq m-1, \;
\lim\limits_{\nu\rightarrow\nu_0^+}\lambda_{n}({\nu})=+\infty,\; m\leq n \leq k-1.
\vspace{-0.5cm}$$
  \item[{\rm (iv)}] If the $n$-th eigenvalue functions $\lambda_n({\nu})$ are non-decreasing
 in  $ (\nu_0,\nu_0+\epsilon_0)$ for all $0\leq n \leq k-1$, then
they
 satisfy the following asymptotic behaviors  near $\nu_0$:
    \vspace{-0.2cm}$$
\lim\limits_{\nu\rightarrow\nu_0^+}\lambda_{n}({\nu})=-\infty,\;\;0\leq n \leq k-m-1,\;
\lim\limits_{\nu\rightarrow\nu_0^+}\lambda_n({\nu})=\lambda_{n-k+m}({\nu_0}), \;\;k-m\leq n\leq k-1.
 \vspace{-0.2cm}$$
\end{itemize}}

\noindent{\bf Proof.} We only show that (i) holds. The other claims  can be shown similarly.

Let $(r_1,r_2)$ be a finite interval such that $\lambda_j(\nu_0)\in(r_1,r_2)$ for
$0\leq j\leq m-1$. Then, by Lemma 2.5 there exists  $0<\epsilon_0<\epsilon$  such that for each
 $\nu\in(\nu_0-\epsilon_0,\nu_0)$,  $(\pmb\omega,\mathbf{A})_{\nu}$  has exactly  $m$ eigenvalues  in $[r_1,r_2]$, which are
all in $(r_1,r_2)$. Its other $k-m$ eigenvalues out of $[r_1,r_2]$ are denoted by $\hat\lambda_1(\nu)\leq\cdots\leq\hat
\lambda_{k-m}(\nu)$.

Fix $1\leq i\leq k-m$. We claim that $\mathcal{O}^-:=\{\hat\lambda_i({\nu}):\nu\in(\nu_0-\epsilon_0,\nu_0)\}\subset(-\infty,r_1)$.
Otherwise,  by  (i) of Theorem 2.3, $\mathcal{O}^-\subset(r_2,+\infty)$
and there exists  $0\leq i_0 \leq k-1$ such that $\hat\lambda_{i}({\nu})=\lambda_{i_0} ({\nu})$ for each
$\nu\in(\nu_0-\epsilon_0,\nu_0)$. Since  $\lambda_{i_0}(\nu)$
 is non-increasing in $(\nu_0-\epsilon_0,\nu_0)$, $\hat\lambda_{i}({\nu})=\lambda_{i_0}({\nu})
\leq\lambda_{i_0}({\nu_0-\epsilon_0/2})=\hat\lambda_{i}({\nu_0-\epsilon_0/2})$ for each $\nu\in(\nu_0
-\epsilon_0/2, \nu_0)$. By Lemma 2.5, there exists $0<\epsilon_1<\epsilon_0/2$ such that $(\pmb\omega,\mathbf{A})_{\nu}$ has exactly $m$ eigenvalues in
$(r_1, \hat\lambda_{i}({\nu_0-\epsilon_0/2})+1)$ for each $\nu\in(\nu_0-\epsilon_1, \nu_0]$.  Since $r_2<\hat\lambda_{i}
({\nu})<\hat\lambda_{i}({\nu_0-\epsilon_0/2})+1$  and $(\pmb\omega,\mathbf{A})_{\nu}$ has exactly
$m$ eigenvalues in $(r_1,r_2)$  for each $\nu\in(\nu_0-\epsilon_1, \nu_0)$, it follows that
$(\pmb\omega,\mathbf{A})_{\nu}$ has at least $m+1$ eigenvalues in $(r_1, \hat\lambda_{i}({\nu_0-\epsilon_0/2})+1)$ for each $\nu\in(\nu_0-\epsilon_1, \nu_0)$. This is a contradiction. Hence the claim holds.

Consequently, $\hat\lambda_n({\nu})=\lambda_{n-1}({\nu})$, $\nu\in(\nu_0-\epsilon_0,\nu_0)$, for each $1\leq n\leq k-m$ again by (i) of Theorem 2.3.
It follows from (ii) of Theorem 2.3 that (2.12) holds and then (2.13) holds by (i) of Theorem 2.2. The proof is complete.\medskip

\noindent{\bf Remark 2.3.} If $m=0$ and the conditions in (i) of Lemma 2.7 are satisfied, then
the $n$-th eigenvalue functions $\lambda_n({\nu})$, $0\leq n \leq k-1$,
 satisfy the following asymptotic behaviors near $\nu_0$:
\vspace{-0.2cm} $$\lim\limits_{\nu\rightarrow\nu_0^-}\lambda_n({\nu})=-\infty,\; 0\leq n \leq k-1.        \vspace{-0.2cm} $$
     \vspace{-0.2cm}
(ii)--(iv) in Lemma 2.7 can be modified similarly in the case that $m=0$.\medskip

\noindent{\bf Remark 2.4.} If the conditions in (i) and (iii) (or (ii) and (iv)) of Lemma 2.7 are satisfied, then $(\pmb\omega,\mathbf{A})_{\nu_0}$ is a discontinuity point of  $\lambda_n$
for all $0\leq n \leq k-1$.\medskip

\noindent{\bf Remark 2.5.} In Example 2.1, $\lambda_n(s)$ is non-decreasing in $s\in[0,1)$, and non-increasing in $s\in(1,2]$ for each $n=0,1$. Thus, the conditions in (ii)--(iii) of Lemma 2.7 are satisfied for $k=2$ and $m=1$, and therefore,
\vspace{-0.2cm}$$\lim\limits_{s\rightarrow1^-}\lambda_0(s)=\lim\limits_{s\rightarrow1^+}\lambda_0(s)=\lambda_0(1),\;
\lim\limits_{s\rightarrow1^-}\lambda_1(s)=+\infty,\;\lim\limits_{s\rightarrow1^+}\lambda_1(s)=+\infty.\vspace{-0.2cm}$$
This also shows that $\lambda_0(s)$ is continuous at $s=1$, and $\lambda_1(s)$ is not continuous at $s=1$.
\bigskip

\noindent{\bf 3. Continuity and discontinuity of the $n$-th eigenvalue function in the space of the SLEs}\medskip

In this section, the continuous and discontinuous dependence of the $n$-th eigenvalue function  on equation (1.1) is discussed.
Its continuity and discontinuity sets in $\Omega_N^{\mathbb R,+}$ are given and its monotonicity in some directions in the continuity set is studied.  In particular, its  asymptotic behaviors  near a discontinuity point are
completely characterized.
\medskip

It was shown in [10]  that the $n$-th eigenvalue depends continuously on the differential equation (1.5) in the continuous case.
However, the following example  shows that the $n$-th eigenvalue may not depend continuously on equation (1.1) in the discrete case.
\medskip

\noindent{\bf Example 3.1.} \rm Let $s\in[1/10,2]$. Consider the 1-parameter family of
the SLPs, in which
\vspace{-0.2cm}$$ f_0={1/ s}, \ f_1=f_2=1, \; q_1 =q_2 =0, \; w_1 =w_2 =1,\; N=2,\vspace{-0.2cm}$$
 and the coefficients in (1.2) is the same as (2.6).
Then, by Lemma 2.3, direct calculations deduce that the
characteristic function is
\vspace{-0.2cm}$$
\Gamma(\lambda) = (s-1)\lambda^2-s\lambda+1-s.
\vspace{-0.2cm}$$
Thus, the  SLP with $s=1$ has exactly one eigenvalue and with each $s\in[1/10,1)\cup(1,2]$ has exactly two eigenvalues, which are given as
\vspace{-0.2cm}$$\begin{array}{cccc}\begin{array}{cccc}
\lambda_0(s) =\begin{cases} {s+\sqrt{5s^2-8s+4} \over 2(s-1) } & \text{ if } s\in [1/10,1), \\
                     0 & \text{ if } s=1, \\
                      {s-\sqrt{5s^2-8s+4} \over 2(s-1)} & \text{ if } s\in (1,2],\end{cases} \end{array}
&\begin{array}{cccc}
\lambda_1(s) =\begin{cases} {s-\sqrt{5s^2-8s+4} \over 2(s-1)} & \text{ if } s\in [1/10,1), \\
                    {s+\sqrt{5s^2-8s+4}) \over 2(s-1)} & \text{ if }  s\in (1,2]. \end{cases}
\end{array}\end{array}\vspace{-0.2cm}$$
See Figure 3.1 below. It is evident that the $n$-th eigenvalue function  $\lambda_n$
is not continuous at $s=1$ for each $n=0,1$.  This example shows that the $n$-th eigenvalue may not depend continuously  on $1/f_0$ in general.
\vspace{-0.2cm}
\begin{center}

\includegraphics[width=80mm]{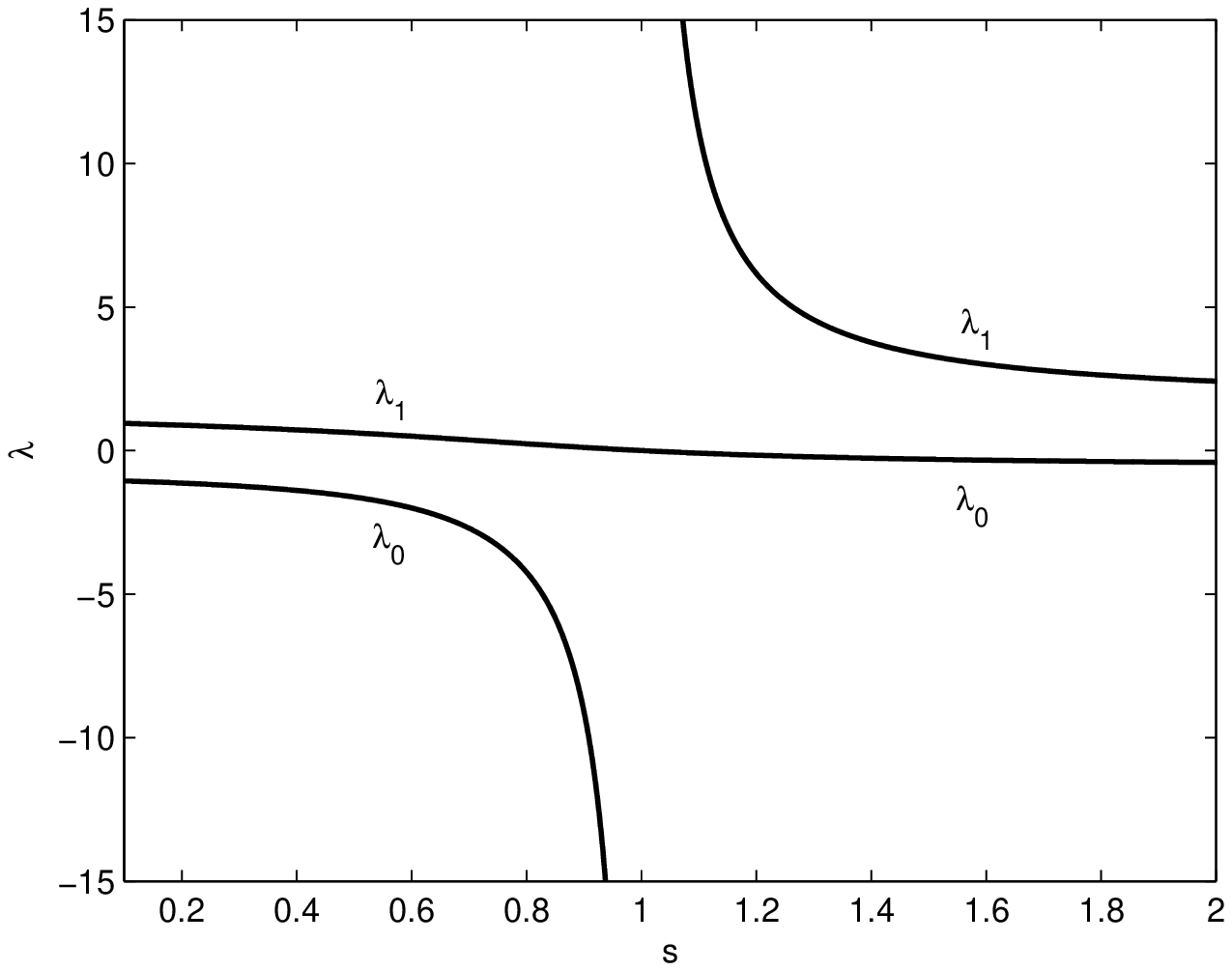}

{\small {\bf Figure 3.1.}}
\end{center}
\vspace{-0.2cm}

Fix a  BC
\vspace{-0.2cm}$$\hat{\mathbf{A}}=\left [\begin{array} {cccc}\hat a_{11}&\hat a_{12}&\hat b_{11}&\hat b_{12}\\
\hat a_{21}&\hat a_{22}&\hat b_{21}&\hat b_{22}\end{array}  \right ]\vspace{-0.1cm}$$
 in this section. We shall study  continuous dependence of
the $n$-th eigenvalue on the SLE.

We need the following monotonicity results of the continuous eigenvalue branches:\medskip

\noindent{\bf Lemma 3.1} {\rm [22, Theorem 4.8]}. {\it Fix a  BC $\hat{\mathbf{A}}$. Then,
each continuous eigenvalue branch $\Lambda$ over
$\Omega_N^{\mathbb R,+}$ is non-increasing in every $(1/f_j)$-direction for
$0 \leq j \leq N-1$, independent of $f_N$, and non-decreasing in every
$q_j$-direction; while its positive and negative  parts are
non-increasing and  non-decreasing in every $w_j$-direction, respectively.}\medskip

For convenience, we introduce the following notations. Let $\mu_1:=\hat a_{11}\hat b_{22}-\hat a_{21}\hat b_{12}$,
 $\mu_2:=\hat a_{22}\hat b_{12}-\hat a_{12}\hat b_{22}$. If $\mu_1\neq 0$, denote
  \vspace{-0.2cm}$$\begin{array}{cccc}\eta:=-\mu_2/\mu_1,\;
   \mathcal{E}:=\{\pmb\omega\in\Omega_N^{\mathbb R,+}: 1/f_0=\eta\}, \\[1.0ex]
   \mathcal{E}^+:=\{\pmb\omega\in\Omega_N^{\mathbb R,+}: 1/f_0>\eta\},\; \mathcal{E}^-:=\{\pmb\omega\in\Omega_N^{\mathbb R,+}: 1/f_0<\eta\}, \\[1.0ex]  \mathcal{E}_1:=\{\pmb\omega\in\Omega_N^{\mathbb R,+}: 1/f_0= -1/\hat a_{11} \} \;{\rm if}\; \hat a_{11}\neq 0;\;\;
\mathcal{E}_2:=\{\pmb\omega\in\Omega_N^{\mathbb R,+}: 1/f_0=\hat a_{12} \}. \end{array}$$
And $\mathcal{E}_1^+$, $\mathcal{E}_1^-$, $\mathcal{E}_2^+$,  $\mathcal{E}_2^-$ can be similarly defined.

Due to (2.5) and the discussion below it, we shall divide our study here into the following three cases:
\vspace{-0.2cm}\begin{itemize}\vspace{-0.2cm}
\item[{\rm (i)}] $\mu_1\neq0,\mu_2\neq0$;\vspace{-0.2cm}
\item[{\rm (ii)}] either $\mu_1=0,\mu_2\neq0$, or $\mu_1\neq0,\mu_2=0$;\vspace{-0.2cm}
\item[{\rm (iii)}] $\mu_1=0,\mu_2=0$.\vspace{-0.2cm}
\end{itemize}

\noindent{\bf Theorem 3.1.} {\it Fix a  BC $\hat{\mathbf{A}}$. Assume that $\mu_1\neq0$ and $\mu_2\neq0$. Then for each $\pmb\omega\in \mathcal{E}$, $(\pmb\omega,\hat{\mathbf{A}})$  has exactly $N-1$ eigenvalues, and for
 each $\pmb\omega\in\Omega_N^{\mathbb R,+}\backslash\mathcal{E}$, $(\pmb\omega,\hat{\mathbf{A}})$  has exactly $N$ eigenvalues $\lambda_n(\pmb\omega):=
 \lambda_n(\pmb\omega,\hat{\mathbf{A}})$, $0\leq n \leq N-1$, which satisfy that
 \begin{itemize}\vspace{-0.2cm}
\item[{\rm (i)}]  $\lambda_n(\pmb\omega)$ are continuous in $\Omega_N^{\mathbb R,+}\backslash\mathcal{E}$;\vspace{-0.2cm}
\item[{\rm (ii)}]  $\lambda_n(\pmb\omega)$ restricted in each connected component of $ \mathcal{E}^+$ or $\mathcal{E}^-$ are non-increasing in every $(1/f_j)$-direction, independent of $f_N$, and non-decreasing in every
$q_j$-direction, while its positive and negative parts are
non-increasing and  non-decreasing in every $w_j$-direction, respectively;\vspace{-0.2cm}\vspace{-0.2cm}
\item[{\rm (iii)}] $\lambda_n(\pmb\omega)$  are not continuous at each point of  $\mathcal{E}$, and have the following asymptotic behaviors
 near any given $\pmb\omega_0\in\mathcal{E}$:
\vspace{-0.2cm}$$\lim\limits_{\mathcal{E}^-\ni\pmb\omega\rightarrow\pmb\omega_0}\lambda_0(\pmb\omega)=-\infty,      \eqno(3.1)\vspace{-0.7cm}$$
\vspace{-0.7cm}$$\lim\limits_{\mathcal{E}^-\ni\pmb\omega\rightarrow\pmb\omega_0}\lambda_n(\pmb\omega)=\lambda_{n-1}(\pmb\omega_0), \;\; 1\leq n \leq N-1,                                                                                                                 \eqno(3.2)\vspace{-0.5cm}$$
\vspace{-0.7cm}$$\lim\limits_{\mathcal{E}^+\cup\mathcal{E}\ni\pmb\omega\rightarrow\pmb\omega_0}\lambda_n(\pmb\omega)=\lambda_{n}(\pmb\omega_0), \;\; 0\leq n \leq N-2,                                                                                                          \eqno(3.3)\vspace{-0.5cm}$$
\vspace{-0.7cm}$$\lim\limits_{\mathcal{E}^+\ni\pmb\omega\rightarrow\pmb\omega_0}\lambda_{N-1}(\pmb\omega)=+\infty.      \eqno(3.4)\vspace{-0.5cm}$$

And consequently, $\lambda_n(\pmb\omega)$ restricted in $\mathcal{E}^+\cup\mathcal{E}$ is continuous  for each $0\leq n \leq N-2$.
\end{itemize}}

\noindent{\bf Proof.}
 By (2.4),  $r(\pmb\omega,\hat{\mathbf{A}})=1$ for each $\pmb\omega\in \mathcal{E}$, and $r(\pmb\omega,\hat{\mathbf{A}})=2$ for each $\pmb\omega\in \Omega_N^{\mathbb R,+}\backslash\mathcal{E}$. It follows from Lemma 2.4 that $(\pmb\omega,\hat{\mathbf{A}})$ has exactly $N-1$
eigenvalues for each $\pmb\omega\in \mathcal{E}$ and $(\pmb\omega,\hat{\mathbf{A}})$ has exactly $N$
eigenvalues for each $\pmb\omega\in\Omega_N^{\mathbb R,+}\backslash \mathcal{E}$.

Now, we show that (i) and (ii) hold.
Since $({\pmb\omega},\hat{\mathbf{A}})$ has exactly $N$
eigenvalues $\lambda_n(\pmb\omega)$, $0\leq n\leq N-1$, for each $\pmb\omega\in \Omega_N^{\mathbb R,+}\backslash\mathcal{E}$, which  is an open set of $\Omega_N^{\mathbb R,+}$,
  $\lambda_n(\pmb\omega)$ is continuous in $\Omega_N^{\mathbb R,+}\backslash\mathcal{E}$  and
 is locally a continuous eigenvalue branch  in each connected component of $\mathcal{E}^+$ or $ \mathcal{E}^-$ for each $0\leq n \leq N-1$ by Theorem 2.1.
This, together with Lemma 3.1, implies its monotonicity   in each connected component of  $ \mathcal{E}^+$
 or $\mathcal{E}^-$.

Next, we show that (iii) holds.
 It suffices to prove (3.1)--(3.4).  Fix any
 \vspace{-0.2cm}$$\pmb\omega_0:=(\eta, 1/f_0^0,\cdots,1/f_N^0,q_1^0,\cdots,q_N^0,w_1^0,\cdots,w_N^0)\in\mathcal{E}.\vspace{-0.2cm}$$
  Let $(r_1,r_2)$ be a finite interval such that $\lambda_j(\pmb\omega_0)\in(r_1,r_2)$, $0\leq j\leq N-2$.
By Lemma 2.5, there exists a neighborhood $\mathcal{U}$  of $\pmb\omega_0$ in $\Omega_N^{\mathbb R,+}$ such that for  each $\pmb\omega\in\mathcal{U}$, $(\pmb\omega,\hat{\mathbf{A}})$ has exactly $N-1$ eigenvalues in $[r_1,r_2]$, which are all in $(r_1,r_2)$. Denote
 \vspace{-0.2cm}$$\mathcal{U}^-:= \mathcal{E}^-\cap \mathcal{U},\; \mathcal{U}^0:= \mathcal{E}\cap \mathcal{U},\; \mathcal{U}^+:= \mathcal{E}^+\cap \mathcal{U}. \vspace{-0.2cm}$$
  Note that
  $\mathcal{U}$ can be chosen sufficiently small such that $\mathcal{U}^-$, $\mathcal{U}^0$ and $\mathcal{U}^+$ are connected.
Since $(\pmb\omega,\hat{\mathbf{A}})$ has exactly $N$ eigenvalues for each  $\pmb\omega\in\mathcal{U}^-$,
 it has exactly one eigenvalue, denoted by $\hat \lambda(\pmb\omega)$,
outside $[r_1,r_2]$.  It follows from (i) of Theorem 2.3 that either  $F:=\{\hat \lambda(\pmb\omega):\pmb\omega\in\mathcal{U}^-\}
 \subset(-\infty,r_1)$ or $F\subset(r_2, +\infty)$.

 Denote $\pmb\omega(s):=(s,1/f_1^0,\cdots,1/f_N^0,q_1^0,\cdots,q_N^0,
w_1^0,\cdots,w_N^0)$, $s\in\mathbb{R}$, and $\mathcal{O}:=\{\pmb\omega(s):s\in(\eta-\epsilon,\eta)\}$ with $0<\epsilon<|\eta|$. Then  $\pmb\omega_0=\pmb\omega(\eta)$ and for each $s\in(\eta-\epsilon,\eta)$, $(\pmb\omega(s),\hat{\mathbf{A}})$ has
exactly $N$ eigenvalues $\lambda_n(\pmb\omega(s))$, $0\leq n\leq N-1$, and  $(\pmb\omega(\eta),\hat{\mathbf{A}})$ has exactly
$N-1$ eigenvalues. By (ii), $\lambda_n(\pmb\omega(s))$ is non-increasing in $(\eta-\epsilon,\eta)$ for each
$0\leq n \leq N-1$. Hence, $\lim_{s\rightarrow\eta^-}\lambda_0(\pmb\omega(s))=-\infty$ by (i) of Lemma 2.7. This implies that there exists an $\pmb\omega_1\in\mathcal{U}^-$
such that ${\lambda}_0({\pmb\omega}_1)<r_1$. Hence, $\hat{\lambda}({\pmb\omega}_1)={\lambda}_0({\pmb\omega}_1)$.  Thus, again by (i) of Theorem 2.3,
 $F=\{\lambda_0(\pmb\omega):\pmb\omega\in\mathcal{U}^-\}\subset(-\infty,r_1)$.
 By  (ii) of Theorem 2.3 one gets that (3.1) holds.

With a similar argument to the proof of (3.1), one can show that   (3.4) holds.

Note that $\mathcal{U}^0$ is connected and
$(\pmb\omega,\hat{\mathbf{A}})$ has exactly $N-1$ eigenvalues for each $\pmb\omega\in\mathcal{U}^0$. Hence,  $\lambda_n$ restricted in $\mathcal{U}^0$  is continuous and locally a continuous eigenvalue branch for each $0\leq n\leq N-2$ by Theorem 2.1. This, together with  Theorem 2.2, implies that (3.2) and (3.3) hold.
The entire proof is complete.
\medskip

\noindent{\bf Theorem 3.2.} {\it Fix a  BC $\hat{\mathbf{A}}$.
Assume that either $\mu_1=0, \mu_2\neq0$ or $\mu_1\neq0, \mu_2=0$. Then for
each  $\pmb\omega\in\Omega_N^{\mathbb R,+} $,  $(\pmb\omega,\hat{\mathbf{A}})$ has exactly $N$ eigenvalues $\lambda_n(\pmb\omega)$, $0\leq n \leq N-1$, which satisfy that
 \begin{itemize}\vspace{-0.2cm}
\item[{\rm (i)}] $\lambda_n(\pmb\omega)$  are continuous in $\Omega_N^{\mathbb R,+}$;\vspace{-0.2cm}
\item[{\rm (ii)}] $\lambda_n(\pmb\omega)$ restricted in each connected component of $\Omega_N^{\mathbb R,+}$   have the same monotonicity as that in  {\rm(ii)} of Theorem {\rm 3.1}.
 \end{itemize}}

\noindent{\bf Theorem 3.3.} {\it Fix a  BC $\hat{\mathbf{A}}$. Assume that $\mu_1=0$ and $\mu_2=0$. Then the BC $\hat{\mathbf{A}}$ can be written as
\vspace{-0.2cm}$$\vspace{-0.2cm}\begin{array} {llll} either &
 \hat{\mathbf{A}}_1=\left [\begin{array} {llll}\hat a_{11}&-1&0&0\\
0&0&-1&0\end{array}  \right ]  &  or &
 \hat{\mathbf{A}}_2=\left [\begin{array} {llll}1&\hat a_{12}&0&0\\
0&0&-1&0\end{array}  \right ].
\end{array}                                    \vspace{-0.2cm}$$\vspace{-0.2cm}
\begin{itemize}\vspace{-0.2cm}
 \item[{\rm (i)}]\vspace{-0.2cm}
    In the case that  $\hat{\mathbf{A}}=\hat{\mathbf{A}}_1$ with $\hat a_{11}\neq 0$, $(\pmb\omega,\hat{\mathbf{A}})$  has exactly $N-2$ eigenvalues for  each  $\pmb\omega\in \mathcal{E}_1$, and
 $(\pmb\omega,\hat{\mathbf{A}})$  has exactly $N-1$ eigenvalues
 $\lambda_n(\pmb\omega)$, $0\leq n \leq N-2$, for each $\pmb\omega\in \Omega_N^{\mathbb R,+}\backslash\mathcal{E}_1$,  which satisfy that
  \begin{itemize}\vspace{-0.2cm}
\item[{\rm (ia)}] $\lambda_n(\pmb\omega)$  are continuous in $\Omega_N^{\mathbb R,+}\backslash\mathcal{E}_1$;\vspace{-0.2cm}
\item[{\rm (ib)}] $\lambda_n(\pmb\omega)$ restricted in each connected component of  $ \mathcal{E}_1^+$ or $\mathcal{E}_1^-$ have the same monotonicity as that in {\rm(ii)} of Theorem {\rm 3.1};\vspace{-0.2cm}
\item[{\rm (ic)}] $\lambda_n(\pmb\omega)$ are not continuous at each point of  $\mathcal{E}_1$, and have the following asymptotic behaviors
near any given $\pmb\omega_0\in\mathcal{E}_1$:
\vspace{-0.2cm}$$\begin{array} {llll}\lim\limits_{\mathcal{E}_1^-\ni\pmb\omega\rightarrow\pmb\omega_0}\lambda_0(\pmb\omega)=-\infty,
\lim\limits_{\mathcal{E}_1^-\ni\pmb\omega\rightarrow\pmb\omega_0}\lambda_n(\pmb\omega)=\lambda_{n-1}(\pmb\omega_0),\; 1\leq n \leq N-2,\\[1.0ex]
\lim\limits_{\mathcal{E}_1^+\cup\mathcal{E}_1\ni\pmb\omega\rightarrow\pmb\omega_0}\lambda_n(\pmb\omega)=\lambda_{n}(\pmb\omega_0),\;  0\leq n \leq N-3,
\lim\limits_{\mathcal{E}_1^+\ni\pmb\omega\rightarrow\pmb\omega_0}\lambda_{N-2}(\pmb\omega)=+\infty.
\end{array}                                                                                     \vspace{-0.2cm}$$
And consequently, $\lambda_n(\pmb\omega)$ restricted in $\mathcal{E}_1^+\cup\mathcal{E}_1$ is continuous  for each $0\leq n \leq N-3$.
 \end{itemize}
\vspace{-0.2cm}
\item[{\rm (ii)}]
In the case that  $\hat{\mathbf{A}}=\hat{\mathbf{A}}_2$ with $\hat a_{12}\neq 0$,
similar results in {\rm(i)} hold for $\mathcal{E}_1$, $\mathcal{E}_1^+$, $\mathcal{E}_1^-$ replaced by $\mathcal{E}_2$, $\mathcal{E}_2^+$, $\mathcal{E}_2^-$, respectively.\vspace{-0.2cm}
  \item [{\rm (iii)}]
 In the case that $\hat{\mathbf{A}}=\hat{\mathbf{A}}_1$ with $\hat a_{11}=0$ or $\hat{\mathbf{A}}=\hat{\mathbf{A}}_2$ with $\hat a_{12}=0$,
 $(\pmb\omega,\hat{\mathbf{A}})$  has exactly $N-1$ eigenvalues $\lambda_n(\pmb\omega)$, $0\leq n \leq N-2$, for each $\pmb\omega\in \Omega_N^{\mathbb R,+}$,
 which satisfy that
  \begin{itemize}\vspace{-0.2cm}
\item[{\rm (iiia)}] $\lambda_n(\pmb\omega)$ are continuous  in $\Omega_N^{\mathbb R,+}$;\vspace{-0.2cm}
\item[{\rm (iiib)}] $\lambda_n(\pmb\omega)$ restricted in each connected component of $\Omega_N^{\mathbb R,+}$ have the same monotonicity  as that in {\rm (ii)} of Theorem {\rm 3.1}.\vspace{-0.2cm}
\end{itemize}
\end{itemize}}
Since the proofs of Theorems 3.2--3.3 are similar to that of Theorem 3.1, we omit their details.
\bigskip

\noindent{\bf 4. Continuity and discontinuity of the $n$-th eigenvalue function in the space of the BCs}\medskip

In this section, the continuous and discontinuous dependence of the $n$-th eigenvalue function on the boundary condition (1.2) is investigated.
Its continuity and discontinuity sets  in $ \mathcal B^{\mathbb C}$ are given and its monotonicity in some directions in the continuity set is studied.
Especially, its asymptotic behaviors near a discontinuity point are completely characterized.

Fix a difference equation $\hat{\pmb\omega}=(1/\hat f,\hat q,\hat w)\in \Omega_N^{\mathbb{R,+}}$ in this section. By Lemma 2.1, $\mathcal B^{\mathbb C}$
equals the union of
$\mathcal{O}_{1,3}^{\mathbb{C}}$, $\mathcal{O}_{2,3}^{\mathbb{C}}$, $\mathcal{O}_{1,4}^{\mathbb{C}}$ and $\mathcal{O}_{2,4}^{\mathbb{C}}$, which are four open sets of $\mathcal B^{\mathbb C}$.
Thus, we shall  consider the $n$-th eigenvalue function $\lambda_n$ in
 $\mathcal{O}_{1,3}^{\mathbb{C}}$,  $\mathcal{O}_{2,3}^{\mathbb{C}}$, $\mathcal{O}_{1,4}^{\mathbb{C}}$, and  $\mathcal{O}_{2,4}^{\mathbb{C}}$, separately. We shall remark that the method used here is different from that used in the continuous case [10], where the authors divided the BCs  into  the separated and coupled ones. The method used here is more convenient in
dealing with the discrete case. Finally, we shall apply our results to
 the separated and  coupled BCs.

We now introduce the following notations for convenience:
\vspace{-0.05cm}$$\begin{array} {cccc}\mathcal{B}_{1,4}:=\left\{ \mathbf{A}\in\mathcal{O}_{1,4}^{\mathbb{C}}:  a_{12}=1/\hat f_0\right\},\;
\mathcal{B}_{2,4}:=\left\{ \mathbf{A}\in\mathcal{O}_{2,4}^{\mathbb{C}}:  a_{11}=-\hat f_0\right\},\\[2.0ex]
\mathcal{B}_{1,3}:=\left\{ \mathbf{A}\in\mathcal{O}_{1,3}^{\mathbb{C}}:  (a_{12}-1/\hat f_0)b_{22}=|z|^2\right\},\;
\mathcal{B}_{2,3}:=\left\{ \mathbf{A}\in\mathcal{O}_{2,3}^{\mathbb{C}}:  (a_{11}+\hat f_0)b_{22}=|z|^2\right\},\\[2.0ex]
\mathbf{C}:=\left [\begin{array} {cccc}1&1/\hat f_0&0&0\\
0&0&1&0\end{array}  \right ],\;
\mathcal{B}_{1,3r}:=\left\{ \mathbf{A}\in\mathcal{B}_{1,3}: a_{12}\geq 1/\hat f_0,b_{22}\geq 0\right\}\backslash\left\{\mathbf{C}\right\},\\[2.0ex]
\mathcal{B}_{1,3l}:=\left\{ \mathbf{A}\in\mathcal{B}_{1,3}:  a_{12}\leq 1/\hat f_0,b_{22}\leq 0\right\}\backslash\left\{\mathbf{C}\right\},\end{array}\vspace{-0.2cm}$$
\vspace{-0.1cm}$$\begin{array} {cccc}
\mathcal{B}_{1,3r}^+:=\left\{ \mathbf{A}\in\mathcal{O}_{1,3}^{\mathbb{C}}: a_{12}\geq 1/\hat f_0,b_{22}\geq 0, (a_{12}-1/\hat f_0)b_{22}>|z|^2\right\},\\[2.0ex]
\mathcal{B}_{1,3l}^+:=\left\{ \mathbf{A}\in\mathcal{O}_{1,3}^{\mathbb{C}}:  a_{12}\leq 1/\hat f_0,b_{22}\leq 0,  (a_{12}-1/\hat f_0)b_{22}>|z|^2\right\},\\[2.0ex]
\mathcal{B}_{2,3r}:=\left\{ \mathbf{A}\in\mathcal{B}_{2,3}: a_{11}+ \hat f_0\geq 0, b_{22}\geq 0\right\}\backslash\left\{\mathbf{C}\right\},\\[2.0ex]
\mathcal{B}_{2,3l}:=\left\{ \mathbf{A}\in\mathcal{B}_{2,3}:  a_{11}+ \hat f_0\leq 0, b_{22}\leq 0\right\}
\backslash\left\{\mathbf{C}\right\},
\end{array}\vspace{-0.1cm}$$
$\mathcal{B}_{1,4}^+$ and $\mathcal{B}_{2,4}^+$, and
$\mathcal{B}_{1,4}^-$, $\mathcal{B}_{2,4}^-$, $\mathcal{B}_{1,3}^-$ and $\mathcal{B}_{2,3}^-$  can be defined similarly as  $\mathcal{E}^+$ and $\mathcal{E}^-$ in Section 3, respectively;
 $\mathcal{B}_{2,3r}^+$ and $\mathcal{B}_{2,3l}^+$ can be defined similarly as $\mathcal{B}_{1,3r}^+$ and $\mathcal{B}_{1,3l}^+$, respectively.
Then $\mathcal{O}_{i,4}^{\mathbb{C}}=\mathcal{B}_{i,4}^+\cup\mathcal{B}_{i,4}\cup\mathcal{B}_{i,4}^-$, $\mathcal{O}_{i,3}^{\mathbb{C}}=
\mathcal{B}_{i,3r}^+\cup\mathcal{B}_{i,3}\cup\mathcal{B}_{i,3}^-\cup\mathcal{B}_{i,3l}^+$, $\mathcal{B}_{i,3}=\mathcal{B}_{i,3r}\cup\{\mathbf{C}\}\cup\mathcal{B}_{i,3l}$ and $\mathcal{B}_{i,3r}\cap\mathcal{B}_{i,3l}=\varnothing$ for $i=1,2$.

Let $\mathcal{B}:=\cup_{i=1}^2\cup_{j=3}^4\mathcal{B}_{i,j}$. Then  $\theta(\hat{\pmb\omega},\mathbf{A})=0$ if and only if $\mathbf{A}\in\mathcal{B}$ for the fixed equation $\hat{\pmb\omega}$, where $\theta$ is defined by (2.5).\medskip

 Note that there are two real parameters in each $\mathcal{O}_{i,j}^{\mathbb{C}}$, $1\leq i\leq2$, $3\leq j\leq 4$. Now we recall the
 monotonicity  of the continuous eigenvalue branches with respect to the two real parameters, which was obtained in [22]. \medskip

\noindent{\bf Lemma 4.1} {\rm [22, Theorem 4.6]}. {\it Fix a difference equation $\hat{\pmb\omega}$.
 Then, in each of the coordinate systems
$\mathcal O_{1,3}^{\mathbb C}$, $\mathcal O_{1,4}^{\mathbb C}$, $\mathcal
O_{2,3}^{\mathbb C}$ and $\mathcal O_{2,4}^{\mathbb C}$ in $\mathcal B^{\mathbb C}$,
every continuous eigenvalue branch is always non-decreasing in the two
real axis directions}.\medskip

For example, in $\mathcal O_{1,4}^{\mathbb C}$, every continuous eigenvalue
branch is always non-decreasing in the $a_{12}$-direction and in the
$b_{21}$-direction.\medskip

\noindent{\bf Theorem 4.1.} {\it Fix a difference equation $\hat{\pmb\omega}$. Then  for each  $\mathbf{A}\in\mathcal{B}_{1,4}$, $(\hat{\pmb\omega}, \mathbf{A})$ has exactly $N-1$ eigenvalues and for each $\mathbf{A}\in\mathcal O_{1,4}^{\mathbb C}\backslash\mathcal{B}_{1,4}$, $(\hat{\pmb\omega}, \mathbf{A})$
 has exactly $N$ eigenvalues $\lambda_n(\mathbf{A}):=\lambda_n(\hat{\pmb\omega},\mathbf{A})$, $0\leq n \leq N-1$, which satisfy that
\begin{itemize}\vspace{-0.3cm}
\item[{\rm (i)}] $\lambda_n(\mathbf{A})$   are continuous in $ \mathcal O_{1,4}^{\mathbb C}\backslash\mathcal{B}_{1,4}$;\vspace{-0.3cm}
\item[{\rm (ii)}] $\lambda_n(\mathbf{A})$ restricted in $\mathcal{B}_{1,4}^+$ and $\mathcal{B}_{1,4}^-$ are always non-decreasing in the $a_{12}$-direction and in the
$b_{21}$-direction;\vspace{-0.3cm}
 \item[{\rm (iii)}] $\lambda_n(\mathbf{A})$ are not continuous at each point of $\mathcal{B}_{1,4}$ and  have the following asymptotic behaviors  near any given $\mathbf{A}_0\in\mathcal{B}_{1,4}$:\vspace{-0.2cm}
$$\vspace{-0.5cm}
\lim\limits_{\mathcal{B}_{1,4}^-\cup\mathcal{B}_{1,4}\ni\mathbf{A}\rightarrow \mathbf{A}_0}\lambda_{n}(\mathbf{A})=\lambda_{n}(\mathbf{A}_0),\;0\leq n \leq N-2,\;\lim\limits_{\mathcal{B}_{1,4}^-\ni\mathbf{A}\rightarrow \mathbf{A}_0}\lambda_{N-1}(\mathbf{A})=+\infty,\eqno(4.1)$$\vspace{-0.2cm}
$$\vspace{-0.2cm}\lim\limits_{\mathcal{B}_{1,4}^+\ni\mathbf{A}\rightarrow \mathbf{A}_0}\lambda_{0}(\mathbf{A})=-\infty,\;\lim\limits_{\mathcal{B}_{1,4}^+\ni\mathbf{A}\rightarrow \mathbf{A}_0}\lambda_{n}(\mathbf{A})=\lambda_{n-1}(\mathbf{A}_0),\;\;1\leq n \leq N-1.
                                                                                                          \eqno(4.2)\vspace{-0.5cm}$$
Consequently,   $\lambda_n(\mathbf{A})$ restricted in $\mathcal{B}_{1,4}^-\cup\mathcal{B}_{1,4}$
is  continuous  for each $0\leq n \leq N-2$.\vspace{-0.2cm}
\end{itemize}}

\noindent{\bf Proof.} By (2.4),  $r(\hat{\pmb\omega},\mathbf{A})=1$ for each $\mathbf{A}\in \mathcal{B}_{1,4}$, and $r(\hat{\pmb\omega},\mathbf{A})=2$ for each $\mathbf{A}\in\mathcal O_{1,4}^{\mathbb C}\backslash \mathcal{B}_{1,4}$. It follows from Lemma 2.4 that $(\hat{\pmb\omega},\mathbf{A})$ has exactly $N-1$
eigenvalues for each $\mathbf{A}\in \mathcal{B}_{1,4}$ and $(\hat{\pmb\omega},\mathbf{A})$ has exactly $N$
eigenvalues for each $\mathbf{A}\in \mathcal O_{1,4}^{\mathbb C}\backslash \mathcal{B}_{1,4}$.

Now, we show that (i) and (ii) hold. Since $(\hat{\pmb\omega},\mathbf{A})$ has exactly $N$
eigenvalues $\lambda_n(\mathbf{A})$, $0\leq n\leq N-1$, for each $\mathbf{A}\in \mathcal{B}_{1,4}^+$ (or $\mathcal{B}_{1,4}^-$), which  is an open and connected subset of $\mathcal B^{\mathbb C}$,
 $\lambda_n(\mathbf{A})$  is continuous  and
  locally a continuous eigenvalue branch  in $\mathcal{B}_{1,4}^+$ (or $ \mathcal{B}_{1,4}^-$) for each $0\leq n \leq N-1$  by Theorem 2.1.
This, together with Lemma 4.1, implies its monotonicity   in  $\mathcal{B}_{1,4}^+$ (or $ \mathcal{B}_{1,4}^-$).

Next, we show that (iii) holds. It suffices to  show that (4.1)--(4.2) hold for any given $\mathbf{A}_0\in\mathcal{B}_{1,4}$.
 Let $(r_1,r_2)$ be a finite interval such that $\lambda_j(\mathbf{A}_0)\in (r_1,r_2)$ for all $0\leq j \leq N-2$.
 By Lemma 2.5, there exists a neighborhood $\mathcal{V}$ of
 $\mathbf{A}_0$  in $\mathcal{O}_{1,4}^{\mathbb{C}}$ such that for each  $\mathbf{A}\in\mathcal{V}$, $(\hat{\pmb\omega},\mathbf{A})$ has exactly $N-1$ eigenvalues in $[r_1,r_2]$ that are all in $(r_1,r_2)$.
  Let
  \vspace{-0.2cm}$$\mathcal{V}^-:= \mathcal{B}_{1,4}^-\cap \mathcal{V},\;\mathcal{V}^0:= \mathcal{B}_{1,4}\cap \mathcal{V},\;\mathcal{V}^+:= \mathcal{B}_{1,4}^+\cap \mathcal{V}.\vspace{-0.2cm}$$
Note that $\mathcal{V}$ can be chosen such that $\mathcal{V}^-$, $\mathcal{V}^0$ and $\mathcal{V}^+$ are connected.
   Since  $(\hat{\pmb\omega},\mathbf{A})$  has exactly $N$ eigenvalues for each $\mathbf{A}\in\mathcal{V}^-$,
it has exactly an eigenvalue, denoted by $\hat \lambda(\mathbf{A})$,  outside $[r_1,r_2]$.
  By (i) of Theorem 2.3, either $G:=\{\hat \lambda(\mathbf{A}):\mathbf{A}\in\mathcal{V}^-\}\subset(-\infty,r_1)$ or $G\subset(r_2, +\infty)$.

Suppose that \vspace{-0.2cm}$$\mathbf{A}_0=\left [\begin{array} {cccc}1&1/\hat f_0&\bar{z}^0&0\\
0&z^0&b_{21}^0&1\end{array}  \right ]\in\mathcal{B}_{1,4}.\vspace{-0.2cm}$$
 Denote
 \vspace{-0.2cm}$$\mathbf{A}(s):=\left [\begin{array} {cccc}1&s&\bar{z}^0&0\\
0&z^0&b_{21}^0&1\end{array}  \right ],\; s\in\mathbb{R}. \vspace{-0.2cm}$$
Then  $\mathbf{A}_0=\mathbf{A}(1/\hat f_0)$.
Note that each $(\hat{\pmb\omega},\mathbf{A}(s))$ has exactly $N$ eigenvalues for each $s\in(-\infty,1/\hat f_0)$, and
$(\hat{\pmb\omega},\mathbf{A}(1/\hat f_0))$ has exactly $N-1$ eigenvalues.
 By (ii), $\lambda_n(\mathbf{A}(s))$ is non-decreasing in $(-\infty,1/\hat f_0)$ for each
$0\leq n \leq N-1$.
Hence, $\lim_{s\rightarrow1/\hat f_0^-}\lambda_{N-1}({\mathbf{A}}(s))=+\infty$ by (ii) of Lemma 2.7. This implies that there exists an $\mathbf{A}_1\in\mathcal{V}^-$
such that ${\lambda}_{N-1}({\mathbf{A}}_1)>r_2$. Hence, $\hat{\lambda}({\mathbf{A}}_1)={\lambda}_{N-1}({\mathbf{A}}_1)$.  Thus, again by (i) of Theorem 2.3,
 $G=\{\lambda_{N-1}(\mathbf{A}):\mathbf{A}\in\mathcal{V}^-\}\subset(r_2,+\infty)$.
By  (ii) of Theorem 2.3, the second relation in (4.1) holds.

With a similar argument to the proof of  the second relation in (4.1), one can show that the first relation in (4.2) holds.

Since $\mathcal{V}^0$ is connected and
$(\hat{\pmb\omega},\mathbf{A})$ has exactly $N-1$ eigenvalues for each  $\mathbf{A}\in\mathcal{V}^0$,
by Theorem 2.1 $\lambda_n$ restricted in $\mathcal{V}^0$ is continuous and locally a continuous eigenvalue branch for each $0\leq n\leq N-2$. This, together with Theorem 2.2, shows that the first relation in (4.1) and the second relation in (4.2) hold. This completes the proof.
\medskip

\noindent{\bf Theorem 4.2.} {\it  Fix a difference equation $\hat{\pmb\omega}$.
Then similar results in Theorem {\rm 4.1} hold, where  $\mathcal{O}^{\mathbb{C}}_{1,4}$, $\mathcal{B}_{1,4}$, $\mathcal{B}_{1,4}^+$, $\mathcal{B}_{1,4}^-$, and $a_{12}$ are replaced by
$\mathcal{O}^{\mathbb{C}}_{2,4}$, $\mathcal{B}_{2,4}$, $\mathcal{B}_{2,4}^+$, $\mathcal{B}_{2,4}^-$, and $a_{11}$, respectively.
The corresponding relations in {\rm(4.1)--(4.2)} are denoted by {\rm(4.1$'$)--(4.2$'$)}.}\medskip

\noindent{\bf Proof.} Since the proof  is similar to that of Theorem 4.1, we omit its details.\medskip

\noindent{\bf Theorem 4.3.} {\it Fix a difference equation $\hat{\pmb\omega}$. Then   $(\hat{\pmb\omega}, \mathbf{C})$ has exactly $N-2$ eigenvalues, and for each  $\mathbf{A}\in\mathcal{B}_{1,3}\backslash\{\mathbf{C}\}$, $(\hat{\pmb\omega}, \mathbf{A})$ has exactly $N-1$ eigenvalues, and  for each $\mathbf{A}\in\mathcal O_{1,3}^{\mathbb C}\backslash\mathcal{B}_{1,3}$,
 $(\hat{\pmb\omega}, \mathbf{A})$
 has exactly $N$ eigenvalues $\lambda_n(\mathbf{A})$, $0\leq n \leq N-1$, which satisfy that
\begin{itemize}\vspace{-0.2cm}
\item[{\rm (i)}] $\lambda_n(\mathbf{A})$  are continuous in $ \mathcal O_{1,3}^{\mathbb C}\backslash\mathcal{B}_{1,3}$;\vspace{-0.3cm}
\item[{\rm (ii)}] $\lambda_n(\mathbf{A})$ restricted in $\mathcal{B}_{1,3}^-$, $\mathcal{B}_{1,3r}^+$ and $\mathcal{B}_{1,3l}^+$ are non-decreasing
in the $a_{12}$-direction and in the $b_{22}$-direction;\vspace{-0.3cm}
\item[{\rm (iii)}]  $\lambda_n(\mathbf{A})$ are not continuous at each point of $\mathcal{B}_{1,3}$ and furthermore,
\begin{itemize}\vspace{-0.2cm}
\item[{\rm (iiia)}]   they have the following asymptotic behaviors  near  any given $\mathbf{A}_0\in \mathcal{B}_{1,3r}$:
\vspace{-0.2cm}$$\vspace{-0.2cm}\begin{array} {llll}
\lim\limits_{\mathcal{B}_{1,3}^-\cup\mathcal{B}_{1,3r}\ni\mathbf{A}\rightarrow \mathbf{A}_0}\lambda_{n}(\mathbf{A})=\lambda_{n}(\mathbf{A}_0), 0\leq n \leq N-2,\\[1.0ex] \lim\limits_{\mathcal{B}_{1,3}^-\ni\mathbf{A}\rightarrow \mathbf{A}_0}\lambda_{N-1}(\mathbf{A})=+\infty,
\lim\limits_{\mathcal{B}_{1,3r}^+\ni\mathbf{A}\rightarrow \mathbf{A}_0}\lambda_{0}(\mathbf{A})=-\infty,\\[1.0ex]
\lim\limits_{\mathcal{B}_{1,3r}^+\ni\mathbf{A}\rightarrow \mathbf{A}_0}\lambda_{n}(\mathbf{A})=\lambda_{n-1}(\mathbf{A}_0),\;1\leq n \leq N-1;
                                                                                  \end{array}\eqno(4.3)$$
\item[{\rm (iiib)}]  they have the following asymptotic behaviors  near  any given $\mathbf{A}_0\in \mathcal{B}_{1,3l}$:
\vspace{-0.2cm}$$\vspace{-0.2cm}\begin{array} {llll}
\lim\limits_{\mathcal{B}_{1,3}^-\ni\mathbf{A}\rightarrow \mathbf{A}_0}\lambda_{0}(\mathbf{A})=-\infty,\\[1.0ex]
\lim\limits_{\mathcal{B}_{1,3}^-\ni\mathbf{A}\rightarrow \mathbf{A}_0}\lambda_{n}(\mathbf{A})=\lambda_{n-1}(\mathbf{A}_0),\;1\leq n \leq N-1,\\[1.0ex]
\lim\limits_{\mathcal{B}_{1,3l}^+\cup\mathcal{B}_{1,3l}\ni\mathbf{A}\rightarrow \mathbf{A}_0}\lambda_{n}(\mathbf{A})=\lambda_{n}(\mathbf{A}_0),\;\;0\leq n \leq N-2,\\[1.0ex]
\lim\limits_{\mathcal{B}_{1,3l}^+\ni\mathbf{A}\rightarrow \mathbf{A}_0}\lambda_{N-1}(\mathbf{A})=+\infty;
                                                                                      \end{array}\eqno(4.4)$$
\item[{\rm (iiic)}] they have the following asymptotic behaviors  near $\mathbf{C}$:
\vspace{-0.2cm}$$ \vspace{-0.2cm}\begin{array} {llll}\lim\limits_{\mathcal{B}_{1,3r}^+\cup\mathcal{B}_{1,3r}\cup\mathcal{B}_{1,3}^-\ni\mathbf{A}\rightarrow \mathbf{C}}\lambda_{0}(\mathbf{A})=-\infty,\;
\lim\limits_{\mathcal{B}_{1,3r}^+\ni\mathbf{A}\rightarrow \mathbf{C}}\lambda_{1}(\mathbf{A})=-\infty,\\[1.0ex]
\lim\limits_{\mathcal{B}_{1,3l}^+\cup\mathcal{B}_{1,3l}\ni\mathbf{A}\rightarrow \mathbf{C}}\lambda_{n}(\mathbf{A})=\lambda_{n}(\mathbf{C}),\;0\leq n \leq N-3,
\\[1.0ex]
\lim\limits_{\mathcal{B}_{1,3r}\cup\mathcal{B}_{1,3}^-\ni\mathbf{A}\rightarrow \mathbf{C}}\lambda_{n}(\mathbf{A})=\lambda_{n-1}(\mathbf{C}),\;1\leq n \leq N-2,
\\[1.0ex]\lim\limits_{\mathcal{B}_{1,3r}^+\ni\mathbf{A}\rightarrow \mathbf{C}}\lambda_{n}(\mathbf{A})=\lambda_{n-2}(\mathbf{C}),\; 2\leq n \leq N-1,\\[1.0ex]
\lim\limits_{\mathcal{B}_{1,3l}^+\cup\mathcal{B}_{1,3l}\ni\mathbf{A}\rightarrow \mathbf{C}}\lambda_{N-2}(\mathbf{A})=+\infty,\;\lim\limits_{\mathcal{B}_{1,3l}^+\cup\mathcal{B}_{1,3}^-\ni\mathbf{A}\rightarrow \mathbf{C}}\lambda_{N-1}(\mathbf{A})=+\infty.
\end{array}\eqno(4.5)\vspace{-0.2cm}$$
\vspace{-0.2cm}
\end{itemize}\vspace{-0.2cm}
\end{itemize}
And consequently,  $\lambda_n(\mathbf{A})$ restricted  in $\mathcal{B}_{1,3}^-\cup\mathcal{B}_{1,3r}$ and $\mathcal{B}_{1,3l}^+\cup\mathcal{B}_{1,3l}\cup\{\mathbf{C}\}$ is continuous for each $0\leq n \leq N-3$, and $\lambda_{N-2}(\mathbf{A})$ restricted  in $\mathcal{B}_{1,3}^-\cup\mathcal{B}_{1,3r}$ and $\mathcal{B}_{1,3l}^+\cup\mathcal{B}_{1,3l}$ is continuous.}\medskip

\noindent{\bf Proof.} Since the proofs of (i), (ii), (iiia) and (iiib) are similar to those of Theorem 4.1, we omit their details.

The rest is to show that (iiic) holds.
Note that $(\hat{\pmb\omega},\mathbf{C})$ has exactly $N-2$ eigenvalues.
Let $(r_1,r_2)$ be a finite interval such that $\lambda_{j}(\mathbf{C})\in(r_1,r_2)$ for all $0\leq j\leq N-3$.
 By Lemma 2.5, there exists a neighborhood $\mathcal{V}_1$ of $\mathbf{C}$ in $\mathcal{O}_{1,3}^{\mathbb{C}}$  such that for each  $\mathbf{A}\in\mathcal{V}_1$,  $(\hat{\pmb\omega},\mathbf{A})$
 has exactly $N-2$ eigenvalues in $[r_1,r_2]$ that are all in $(r_1,r_2)$.
 Note that $\mathcal{V}_1$ can be chosen such that  $\mathcal{V}_1\cap\mathcal{B}_{1,3r}^+$,
$\mathcal{V}_1\cap\mathcal{B}_{1,3r}$, $\mathcal{V}_1\cap\mathcal{B}_{1,3}^-$, $\mathcal{V}_1\cap\mathcal{B}_{1,3l}$, and
$\mathcal{V}_1\cap\mathcal{B}_{1,3l}^+$ are connected.
 Then we divide our proof in three steps.

{\bf Step 1.} We show that
\vspace{-0.2cm}$$\vspace{-0.5cm}\lim\limits_{\mathcal{B}_{1,3r}\ni\mathbf{A}\rightarrow\mathbf{C}}\lambda_{0}(\mathbf{A})=-\infty,
\lim\limits_{\mathcal{B}_{1,3r}\ni\mathbf{A}\rightarrow \mathbf{C}}\lambda_{n}(\mathbf{A})=\lambda_{n-1}(\mathbf{C}),\;1\leq n \leq N-2,\eqno(4.6)\vspace{-0.2cm}$$\vspace{-0.2cm}
$$\vspace{-0.5cm}\lim\limits_{\mathcal{B}_{1,3l}\ni\mathbf{A}\rightarrow \mathbf{C}}\lambda_{n}(\mathbf{A})=\lambda_{n}(\mathbf{C}),0\leq n \leq N-3, \lim\limits_{\mathcal{B}_{1,3l}\ni\mathbf{A}\rightarrow \mathbf{C}}\lambda_{N-2}(\mathbf{A})=+\infty.\eqno(4.7)$$\vspace{-0.1cm}

Since for each $\mathbf{A}\in\mathcal{V}_1\cap\mathcal{B}_{1,3r}$, $(\hat{\pmb\omega},\mathbf{A})$
has exactly $N-1$ eigenvalues, and then
has exactly an eigenvalue, denoted by $\hat \lambda(\mathbf{A})$,  outside $[r_1,r_2]$.
  By  (i) of Theorem 2.3,  $H:=\{\hat \lambda(\mathbf{A}):\mathbf{A}\in\mathcal{V}_1\cap\mathcal{B}_{1,3r}\}\subset(-\infty,r_1)$ or $H\subset(r_2, +\infty)$.

  Let
  \vspace{-0.2cm}$$\mathbf{A}_1(s)=\left [\begin{array} {cccc}1&s&0&0\\
0&0&-1&0\end{array}  \right ],\;s\in\mathbb{R}.$$
Then $\mathbf{A}_1(1/\hat f_0)=\mathbf{C}$. For each $s\in(1/\hat f_0,+\infty)$, $(\hat{\pmb\omega},\mathbf{A}_1(s))$ has exactly $N-1$ eigenvalues and
$(\hat{\pmb\omega},\mathbf{A}_1(1/\hat f_0))$ has exactly $N-2$ eigenvalues.
 By Theorem 2.1, $\lambda_n(\mathbf{A}_1(s))$  is continuous  and
  locally a continuous eigenvalue branch   in $(1/\hat f_0,+\infty)$ for each $0\leq n \leq N-2$.
This, together with Lemma 4.1, implies that $\lambda_n(\mathbf{A}_1(s))$ is non-decreasing in $(1/\hat f_0,+\infty)$ for each
$0\leq n \leq N-2$.
Hence, $\lim_{s\rightarrow1/\hat f_0^+}\lambda_{0}({\mathbf{A}_1}(s))=-\infty$ by (iv) of Lemma 2.7. This implies that there exists an $\mathbf{A}_1\in\mathcal{V}_1\cap\mathcal{B}_{1,3r}$
such that ${\lambda}_{0}({\mathbf{A}}_1)<r_1$. Hence, $\hat{\lambda}({\mathbf{A}}_1)={\lambda}_{0}({\mathbf{A}}_1)$.  Thus, again by (i) of Theorem 2.3,
 $H=\{\lambda_{0}(\mathbf{A}):\mathbf{A}\in\mathcal{V}_1\cap\mathcal{B}_{1,3r}\}\subset(-\infty,r_1)$.
 Then it follows  from  (ii) of Theorem 2.3 that the first relation in (4.6) holds.

With a similar argument to the proof of the first relation in (4.6), one can show that the second relation in (4.7) holds.

It follows from Theorem 2.2 that the second relation in (4.6) and the first relation in (4.7) hold.

{\bf Step 2.} We show that
\vspace{-0.2cm}$$\vspace{-0.2cm}\lim\limits_{\mathcal{B}_{1,3r}^+\ni\mathbf{A}\rightarrow \mathbf{C}}\lambda_{n}(\mathbf{A})
=-\infty,\;n=0,1,\;\lim\limits_{\mathcal{B}_{1,3r}^+\ni\mathbf{A}\rightarrow \mathbf{C}}\lambda_{n}(\mathbf{A})=\lambda_{n-2}(\mathbf{C}),\;\;2\leq n \leq N-1,\eqno(4.8)\vspace{-0.5cm}$$
  \vspace{-0.5cm}$$\lim\limits_{\mathcal{B}_{1,3l}^+\ni\mathbf{A}\rightarrow \mathbf{C}}\lambda_{n}(\mathbf{A})=\lambda_{n}(\mathbf{C}),0\leq n \leq N-3,\lim\limits_{\mathcal{B}_{1,3l}^+\ni\mathbf{A}\rightarrow\mathbf{C}}\lambda_{n}(\mathbf{A})=+\infty,\;n=N-2,N-1.\eqno(4.9)
\vspace{-0.5cm}$$

Since for each $\mathbf{A}\in\mathcal{V}_1\cap\mathcal{B}_{1,3r}^+$, $(\hat{\pmb\omega},\mathbf{A})$ has exactly $N$ eigenvalues,
 and then has exactly two eigenvalues, denoted by
$\hat \lambda_1(\mathbf{A})\leq\hat \lambda_2(\mathbf{A})$,  outside $[r_1,r_2]$.
By  (i) of  Theorem 2.3, either $I_n:=\{\hat \lambda_n(\mathbf{A}):\mathbf{A}\in\mathcal{V}_1\cap\mathcal{B}_{1,3r}^+\}\subset(-\infty,r_1)$ or $I_n\subset(r_2, +\infty)$ for each $n=1,2$.

Let
\vspace{-0.2cm}$$\mathbf{A}_2(s)=\left [\begin{array} {cccc}1&s+1/\hat f_0&0&0\\
0&0&-1&s\end{array}  \right ],\;s\in\mathbb{R}.$$
 Then $\mathbf{A}_2(0)=\mathbf{C}$.
 $(\hat{\pmb\omega},$ $\mathbf{A}_2(s))$  has exactly $N$ eigenvalues for each $s\in(0,+\infty)$,  and
$(\hat{\pmb\omega},\mathbf{A}_2(0))$ has exactly $N-2$ eigenvalues. By (ii), $\lambda_n(\mathbf{A}_2(s))$ is
non-decreasing in  $(0,+\infty)$ for each $0\leq n \leq N-1$. Hence, by (iv) of Lemma 2.7,
$\lim_{s\rightarrow {0}^+}\lambda_{n}({\mathbf{A}}_2(s))=-\infty,n=0,1.$
This implies that
there exists  $\mathbf{A}_2\in\mathcal{V}_1\cap\mathcal{B}_{1,3r}^+$ such that ${\lambda}_{n}(\mathbf{A}_2)<r_1$, $n=0,1$.
 Hence, $\hat{\lambda}_{n}(\mathbf{A}_2)={\lambda}_{n-1}(\mathbf{A}_2)$, $n=1,2$. Thus, again by (i) of Theorem 2.3,
 $I_n=\{{\lambda}_{n-1}(\mathbf{A}):\mathbf{A}\in\mathcal{V}_1\cap\mathcal{B}_{1,3r}^+\}\subset(-\infty,r_1)$, $n=1,2$.
By  (ii) of Theorem 2.3, the first relation in (4.8) holds.

With a similar argument to the proof of the first relation in (4.8), one can show that the second relation in (4.9) holds.

It follows from Theorem 2.2 that the second relation in  (4.8) and the first relation in (4.9) hold.

{\bf Step 3.} We show that
\vspace{-0.2cm}$$\vspace{-0.2cm}\lim\limits_{\mathcal{B}_{1,3}^-\ni\mathbf{A}\rightarrow \mathbf{C}}\lambda_{0}(\mathbf{A})=-\infty,\;\lim\limits_{\mathcal{B}_{1,3}^-\ni\mathbf{A}\rightarrow \mathbf{C}}\lambda_{N-1}(\mathbf{A})=+\infty,\eqno(4.10)\vspace{-0.7cm}$$
\vspace{-0.7cm}$$\lim\limits_{\mathcal{B}_{1,3}^-\ni\mathbf{A}\rightarrow \mathbf{C}}\lambda_{n}(\mathbf{A})=\lambda_{n-1}(\mathbf{C}),\;\;1\leq N\leq N-2.\eqno(4.11)\vspace{-0.7cm}$$

Since for each $\mathbf{A}\in\mathcal{V}_1\cap\mathcal{B}_{1,3}^-$, $(\hat{\pmb\omega},\mathbf{A})$ has exactly $N$ eigenvalues,
and then has exactly two eigenvalues, denoted by
 $ \lambda'_1(\mathbf{A})\leq \lambda'_2(\mathbf{A})$,  outside $[r_1,r_2]$.
By  (i) of  Theorem 2.3, either $J_n:=\{ \lambda'_n(\mathbf{A}):\mathbf{A}\in\mathcal{V}_1\cap\mathcal{B}_{1,3}^-\}\subset(-\infty,r_1)$ or $J_n\subset(r_2, +\infty)$ for each  $n=1,2$.
Let $\mathbf{A}_3\in \mathcal{V}_1\cap\mathcal{B}_{1,3r}$ and $\mathbf{A}_4\in \mathcal{V}_1\cap\mathcal{B}_{1,3l}$.
Then it follows from Step 1 that  $\lambda_0(\mathbf{A}_3) \in(-\infty,r_1)$ and $\lambda_{N-2}(\mathbf{A}_4)\in(r_2,+\infty)$.
 By Lemma 2.6, $\lambda_0(\mathbf{A}_3)$ lies in a continuous eigenvalue branch $\Lambda_1$ defined in a connected neighborhood $\mathcal{V}_2 $
 of $\mathbf{A}_3$
in $\mathcal{O}_{1,3}^{\mathbb{C}}$.  Thus, one can choose $\mathbf{A}_5\in\mathcal{V}_1\cap\mathcal{V}_2\cap\mathcal{B}_{1,3}^-$
 sufficiently close to $\mathbf{A}_3$
 such that $\Lambda_1(\mathbf{A}_5)\in (-\infty,r_1)$.
Similarly, $\lambda_{N-2}(\mathbf{A}_4)$ lies in a continuous eigenvalue branch $\Lambda_2$ defined in a connected neighborhood $\mathcal{V}_3$
 of
 $\mathbf{A}_4$ in $\mathcal{O}_{1,3}^{\mathbb{C}}$.  Thus, one can choose $\mathbf{A}_6\in\mathcal{V}_1\cap\mathcal{V}_3\cap\mathcal{B}_{1,3}^-$ sufficiently close to
$\mathbf{A}_4$ such that $\Lambda_2(\mathbf{A}_6)\in (r_2,+\infty)$. Since either $J_n\subset(-\infty,r_1)$ or $J_n\subset(r_2, +\infty)$ for each $n=1,2$, again by (i) of Theorem 2.3, $J_1=\{\lambda_0(\mathbf{A}):\mathbf{A}\in\mathcal{V}_1
\cap\mathcal{B}_{1,3}^-\}\subset(-\infty,r_1)$ and
$J_2=\{\lambda_{N-1}(\mathbf{A}):\mathbf{A}\in\mathcal{V}_1\cap\mathcal{B}_{1,3}^-\}\subset(r_2, +\infty)$.
By (ii) of Theorem 2.3,  (4.10) holds.
With a similar argument to that used in the proof of Theorem 2.2, one can show  that (4.11) holds.
The whole proof is complete.
\medskip

\noindent{\bf Theorem 4.4.} {\it Fix a difference equation $\hat{\pmb\omega}$. Then
similar results in Theorem {\rm 4.3} hold, where $\mathcal{O}^{\mathbb{C}}_{1,3}$, $\mathcal{B}_{1,3}$, $\mathcal{B}_{1,3r}$, $\mathcal{B}_{1,3l}$, $\mathcal{B}_{1,3}^-$, $\mathcal{B}_{1,3r}^+$,
$\mathcal{B}_{1,3l}^+$, and $a_{12}$ are replaced by $\mathcal{O}^{\mathbb{C}}_{2,3}$, $\mathcal{B}_{2,3}$, $\mathcal{B}_{2,3r}$, $\mathcal{B}_{2,3l}$, $\mathcal{B}_{2,3}^-$,
$\mathcal{B}_{2,3r}^+$, $\mathcal{B}_{2,3l}^+$, and $a_{11}$, respectively.
The corresponding relations in {\rm(4.3)--(4.5)} are denoted by {\rm(4.3$'$)--(4.5$'$)}.}\medskip

\noindent{\bf Proof.} Since the proof is similar to that of Theorem 4.3, we omit its details.\medskip

Combining Theorems 4.1--4.4 yields the continuity and discontinuity sets for each of the $n$-th eigenvalue function in $\mathcal B^{\mathbb C}$:\medskip

\noindent{\bf Theorem 4.5.} {\it Fix a difference equation $\hat{\pmb\omega}$. Then $(\hat{\pmb\omega}, \mathbf{C})$
has exactly $N-2$ eigenvalues, $(\hat{\pmb\omega}, \mathbf{A})$ has exactly $N-1$ eigenvalues for each  $\mathbf{A}\in\mathcal{B}\backslash\{\mathbf{C}\}$, and $(\hat{\pmb\omega}, \mathbf{A})$
 has exactly $N$ eigenvalues  for each $\mathbf{A}\in\mathcal B^{\mathbb C}\backslash\mathcal{B}$. Moreover, the $n$-th eigenvalue function
$\lambda_n(\mathbf{A})$  is continuous in $ \mathcal B^{\mathbb C}\backslash\mathcal{B}$
 and not continuous at each point of $\mathcal{B}$  for each $0\leq n\leq N-1$. }\medskip

Now we apply Theorems 4.1--4.4 to the separated and coupled BCs, respectively. We  introduce the following notations for convenience:
$\xi:=\arctan(-1/\hat f_0)+\pi$ for $\hat f_0>0$, and $\xi:=\arctan(-1/\hat f_0)$ for $\hat f_0<0$;
\vspace{-0.1cm}$$\begin{array}{cccc} {\mathcal{B}_{S_1}}:=\{\mathbf{S}_{\alpha,\beta}\in \mathcal{B}_S:(\alpha,\beta)\in\{\xi\}\times(0,\pi] \;{\rm or}\; (\alpha,\beta)\in[0,\pi)
\times\{\pi\}\},\\[1.0ex]
{\mathcal{B}_{C_1}}:
=\{[e^{i\gamma}K\,|\,-I]\in\mathcal{B}_{C}:
\gamma\in[0,\pi),\;\;K\in SL(2,\mathbb{R}),\;k_{12}\neq0,\; k_{11}/k_{12}=\hat f_0\},\\[1.0ex]
{\mathcal{B}_{C_1}^+}:=\{[e^{i\gamma}K\,|\,-I]\in\mathcal{B}_{C}:
\gamma\in[0,\pi),\;\;K\in SL(2,\mathbb{R}),\; k_{12}\neq0,\;k_{11}/k_{12}>\hat f_0\},\\[1.0ex]
{\mathcal{B}_{C_1}^-}:=\{[e^{i\gamma}K\,|\,-I]\in\mathcal{B}_{C}:
\gamma\in[0,\pi),\;\;K\in SL(2,\mathbb{R}),\;k_{12}\neq0,\; k_{11}/k_{12}<\hat f_0\},\end{array}\vspace{-0.05cm}$$
where $\mathbf{S}_{\alpha,\beta}$ is defined by (2.2). By Lemma 2.2, $\theta(\hat{\pmb\omega},\mathbf{A})=0$ if and only if $\mathbf{A}\in{\mathcal{B}_{S_1}}\cup{\mathcal{B}_{C_1}}$, where $\theta$ is defined by (2.5).
Thus  \vspace{-0.4cm}$$\mathcal{B}={\mathcal{B}_{S_1}}\cup{\mathcal{B}_{C_1}}.\eqno(4.12)\vspace{-0.9cm}$$

The following result  gives the continuity and discontinuity sets of  $\lambda_n$ restricted in $\mathcal{B}_S$ and  asymptotic behaviors of $\lambda_n$ in $\mathcal{B}^{\mathbb{C}}$  near each discontinuity point.\medskip

\noindent{\bf Corollary 4.1.} {\it Fix a difference equation $\hat{\pmb\omega}$. Then $(\hat{\pmb\omega},\mathbf{S}_{\xi,\pi})$
  has exactly $N-2$ eigenvalues, and $(\hat{\pmb\omega},\mathbf{S}_{\alpha,\beta})$
  has  exactly $N-1$ eigenvalues for each   $\mathbf{S}_{\alpha,\beta}\in\mathcal{B}_{S_1}\backslash\{\mathbf{S}_{\xi,\pi}\}$, and
 $(\hat{\pmb\omega},\mathbf{S}_{\alpha,\beta})$
  has exactly $N$ eigenvalues for each $\mathbf{S}_{\alpha,\beta}\in\mathcal{B}_S\backslash\mathcal{B}_{S_1}$. Moreover,
  the $n$-th eigenvalue function $\lambda_n$  in $\mathcal{B}^{\mathbb{C}}$ is continuous in $\mathcal{B}_S\backslash\mathcal{B}_{S_1}$ and not continuous
 at each point of $\mathcal{B}_{S_1}$ for each $0\leq n\leq N-1$,  and furthermore,\begin{itemize}\vspace{-0.2cm}
\item[{\rm (i)}] for any fixed $\beta_0\in(0,\pi)$, {\rm (4.1)--(4.2)} and {\rm (4.1$'$)--(4.2$'$)} hold for
  $\mathbf{A}_0$ replaced by $\mathbf{S}_{\xi,\beta_0}$;
\item[{\rm (ii)}] for any fixed $\alpha_0\in[0,\pi)\backslash\{\xi,\pi/2\}$,
 {\rm (4.3)} holds for $\mathbf{A}_0$ replaced by $\mathbf{S}_{\alpha_0,\pi}$ in the case that $\mathbf{S}_{\alpha_0,\pi}\in\mathcal{B}_{1,3r}$, and
{\rm (4.4)} holds for $\mathbf{A}_0$ replaced by $\mathbf{S}_{\alpha_0,\pi}$ in the other case that $\mathbf{S}_{\alpha_0,\pi}\in\mathcal{B}_{1,3l}$; \vspace{-0.2cm}
\item[{\rm (iii)}]  {\rm (4.3$'$)} holds for $\mathbf{A}_0$ replaced by $\mathbf{S}_{\pi/2,\pi}$ in the case that $\mathbf{S}_{\pi/2,\pi}\in\mathcal{B}_{2,3r}$, and
{\rm (4.4$'$)} holds for $\mathbf{A}_0$ replaced by $\mathbf{S}_{\pi/2,\pi}$ in the other case that $\mathbf{S}_{\pi/2,\pi}\in\mathcal{B}_{2,3l}$; \vspace{-0.2cm}
\item[{\rm (iv)}]  {\rm (4.5)} and
{\rm (4.5$'$)} hold  for $\mathbf{C}$ replaced by $\mathbf{S}_{\xi,\pi}$.\vspace{-0.2cm}
\end{itemize}}

\noindent{\bf Proof.} The number of eigenvalues of each $(\hat{\pmb\omega},\mathbf{S}_{\alpha,\beta})$ can be easily verified by Lemma 2.4.
Since $\mathcal{B}_S\backslash\mathcal{B}_{S_1}\subset\mathcal{B}^{\mathbb{C}}\backslash\mathcal{B}$, $\lambda_n$  in $\mathcal{B}^\mathbb{C}$ is continuous in $\mathcal{B}_S\backslash\mathcal{B}_{S_1}$ for each $0\leq n\leq N-1$  by Theorem 4.5. Then, it suffices to show that (i)--(iv) hold.

(i) Fix any $\beta_0\in(0,\pi)$. It is clear that
  \vspace{-0.2cm}$$\mathbf{S}_{\xi,\beta_0}=\left [\begin{array} {cccc}1&1/\hat f_0&0&0\\
0&0&-\cot\beta_0&1\end{array}  \right ]\in\mathcal{B}_{1,4}\cap\mathcal{B}_{2,4}.  \vspace{-0.2cm}$$
 Applying   Theorems 4.1-4.2 to $\mathbf{S}_{\xi,\beta_0}$, one gets that  (i) holds.

 (ii) Fix any $\alpha_0\in[0,\pi)\backslash\{\xi,\pi/2\}$. It is clear that
  \vspace{-0.2cm}$$\mathbf{S}_{\alpha_0,\pi}=\left [\begin{array} {cccc}1&-\tan\alpha_0&0&0\\
0&0&-1&0\end{array}  \right ]\in\mathcal{B}_{1,3r}\cup\mathcal{B}_{1,3l}. \vspace{-0.2cm}$$
 Applying (iiia)--(iiib) in Theorem 4.3 to $\mathbf{S}_{\alpha_0,\pi}$,
one gets that  (ii) holds.

(iii)  It is clear that
  \vspace{-0.2cm}$$ \mathbf{S}_{\pi/2,\pi}=\left [\begin{array} {cccc}0&-1&0&0\\
0&0&-1&0\end{array}  \right ]\in\mathcal{B}_{2,3r}\cup\mathcal{B}_{2,3l}.  \vspace{-0.2cm}$$
 Applying  (iiia)--(iiib) in Theorem 4.4 to $\mathbf{S}_{\pi/2,\pi}$,
one gets that  (iii) holds.

(iv)  It is clear that
   \vspace{-0.2cm}$$ \mathbf{S}_{\xi,\pi}=\left [\begin{array} {cccc}1&1/\hat f_0&0&0\\
0&0&-1&0\end{array}  \right ]=\mathbf{C}.  \vspace{-0.2cm}$$
Applying (iiic) in Theorems 4.3--4.4 to $\mathbf{S}_{\xi,\pi}$,
one gets that  (iv) holds.
The proof is complete.\medskip

The following lemma is the monotonicity result of continuous eigenvalue branches with respect to the two real parameters $\alpha$ and $\beta$ for $\mathbf{S}_{\alpha,\beta}$.\medskip

\noindent{\bf Lemma 4.2} {\rm [22, Theorem 4.4]}. {\it   Each continuous eigenvalue branch over
$\mathcal B_{S}$ is always
strictly decreasing in the $\alpha$-direction and always
strictly increasing in the $\beta$-direction}.\medskip

Now, we consider $\lambda_n$ restricted in $\mathcal{B}_{S} $ for each $ 0\leq n\leq N-1$.  For a fixed  $\beta_0\in(0,\pi]$, set $\lambda_n(\alpha):=\lambda_n(\mathbf{S}_{\alpha,\beta_0})$, and for a fixed $\alpha_0\in[0,\pi)$, set
$\lambda_n(\beta):=\lambda_n(\mathbf{S}_{\alpha_0,\beta})$ for convenience.  \medskip
\vspace{-0.5cm}

\noindent{\bf Corollary 4.2.} {\it Fix a difference equation $\hat{\pmb\omega}$. Then\vspace{-0.2cm}
\begin{itemize}\vspace{-0.2cm}
\item[{\rm (i)}] for any fixed $\beta_0\in(0,\pi)$, the $n$-th eigenvalue functions $\lambda_n(\alpha)$  are strictly decreasing in $[0,\xi)$ or $(\xi,\pi)$
 for all $0\leq n\leq N-1$,  and have the following asymptotic behaviors  near $0$ and $\xi$:
 $$\vspace{-0.2cm}\lim\limits_{\alpha\rightarrow\pi^-}\lambda_n(\alpha)=\lambda_n(0),\;\;0\leq n \leq N-1,
                                                                                   \eqno(4.13)\vspace{-0.2cm}$$\vspace{-0.2cm}
$$\vspace{-0.2cm}\begin{array} {llll}
\lim\limits_{\alpha\rightarrow \xi^-}\lambda_{0}(\alpha)=-\infty,\;\;
\lim\limits_{\alpha\rightarrow \xi^-}\lambda_{n}(\alpha)=\lambda_{n-1}(\xi),\;\;1\leq n\leq N-1,\\[2.0ex]
\lim\limits_{\alpha\rightarrow \xi^+}\lambda_{n}(\alpha)=\lambda_{n}(\xi),\;0\leq n \leq N-2,\;\;\lim\limits_{\alpha\rightarrow\xi^+}\lambda_{N-1}(\alpha)=+\infty,                                \end{array}$$
and consequently, $\lambda_n(\alpha)$ is  continuous in $[\xi,\pi)$   for each $0\leq n \leq N-2$;\vspace{-0.2cm}
\item[{\rm (ii)}] for any  fixed $\alpha_0\in[0,\pi)\backslash\{\xi\}$,  the $n$-th eigenvalue functions  $\lambda_n(\beta)$   are strictly increasing in $(0,\pi)$  for all $0\leq n \leq N-1$,
and have the following asymptotic behaviors  near  $\pi$:
$$\vspace{-0.2cm}\begin{array} {llll}\lim\limits_{\beta\rightarrow \pi^-}\lambda_{n}(\beta)=\lambda_{n}(\pi),\;\;0\leq n\leq N-2,\;\;\lim\limits_{\beta\rightarrow \pi^-}\lambda_{N-1}(\beta)=+\infty,\\[1.0ex]
\lim\limits_{\beta\rightarrow 0^+}\lambda_{0}(\beta)=-\infty,\;\;\lim\limits_{\beta\rightarrow 0^+}\lambda_{n}(\beta)=\lambda_{n-1}(\pi),\;\;1\leq n \leq N-1,
\end{array}$$
and consequently, $\lambda_n(\beta)$  is  continuous in $(0,\pi]$   for each $0\leq n\leq N-2$;\vspace{-0.2cm}
\item[{\rm (iii)}]   for $\beta_0=\pi$,  similar results in {\rm(i)} hold for $N-2$, $N-1$ replaced by $N-3$, $N-2$, respectively.
\item[{\rm (iv)}] for $\alpha_0=\xi$, similar results in  {\rm(ii)} hold for $N-2$, $N-1$ replaced by $N-3$, $N-2$, respectively.
\end{itemize}}
\noindent{\bf Proof.} It suffices to show that (4.13) holds since the rest is direct consequence of Theorem 2.1,
 Lemma 2.7, Corollary 4.1 and Lemma 4.2.

Fixed $\beta_0\in(0,\pi)$. Since $\{\mathbf{S}_{\alpha,\beta_0}:\alpha\in[0,\xi)\cup(\xi,\pi)\}$ is connected and  $(\hat{\pmb\omega}, \mathbf{S}_{\alpha,\beta_0})$ has exactly $N$ eigenvalues for each $\alpha\in[0,\xi)\cup(\xi,\pi)$, the $n$-th eigenvalue function $\lambda_n$ is continuous
and locally a continuous eigenvalue branch in $\{\mathbf{S}_{\alpha,\beta_0}:\alpha\in[0,\xi)\cup(\xi,\pi)\}$ for each $0\leq n\leq N-1$  by Theorem 2.1. Then $\lim_{\alpha\rightarrow\pi^-}\mathbf{S}_{\alpha,\beta_0}=\mathbf{S}_{0,\beta_0}$ implies
that (4.13) holds. This completes the proof.\medskip

The following result  gives the continuity and discontinuity sets of  $\lambda_n$ restricted in $\mathcal{B}_C$ and  asymptotic behaviors of $\lambda_n$ in $\mathcal{B}^{\mathbb{C}}$   near each discontinuity point.\medskip

\noindent{\bf Corollary 4.3.} {\it Fix a difference equation $\hat{\pmb\omega}$.
Then $(\hat{\pmb\omega},\mathbf{A})$
  has exactly $N-1$ eigenvalues for each $\mathbf{A}\in\mathcal{B}_{C_1}$, and $(\hat{\pmb\omega},\mathbf{A})$
  has exactly $N$ eigenvalues for each $\mathbf{A}\in\mathcal{B}_{C}\backslash\mathcal{B}_{C_1}$. Moreover, the $n$-th eigenvalue function $\lambda_n$ in $\mathcal{B}^{\mathbb{C}}$ is continuous in $\mathcal{B}_C\backslash\mathcal{B}_{C_1}$, and not continuous at each point of $\mathcal{B}_{C_1}$ for each $0\leq n\leq N-1$; and
 for any given  $\mathbf{A}_1\in{\mathcal{B}_{C_1}}$, {\rm(4.1)--(4.2)} and {\rm(4.1$'$)--(4.2$'$)} hold for $\mathbf{A}_0$ replaced by $\mathbf{A}_1$.}\medskip

\noindent{\bf Proof.} The number of eigenvalues of  $(\hat{\pmb\omega},\mathbf{A})$ can be easily verified by Lemma 2.4 for each  $\mathbf{A}\in \mathcal{B}_{C}$. Since $\mathcal{B}_{C}\backslash\mathcal{B}_{C_1}\subset\mathcal{B}^{\mathbb{C}}\backslash\mathcal{B}$, $\lambda_n$ in $\mathcal{B}^{\mathbb{C}}$ is continuous in $\mathcal{B}_{C}\backslash\mathcal{B}_{C_1}$ for each $0\leq n\leq N-1$ by Theorem 4.5.  Let  $\mathbf{A}_1:=[e^{i\gamma}K|-I]\in {\mathcal{B}_{C_1}}$.
Then $k_{11}\neq 0$ and $k_{12}\neq 0$.  By Lemma 3.18 in [13],
 \vspace{-0.2cm}$$\mathbf{A}_1=\left[\begin{array} {llll}
1&k_{12}/k_{11} & -e^{-i\gamma}/k_{11}&0\\
0&-e^{i\gamma}/k_{11}&-k_{21}/k_{11}&1\end{array}\right]=\left[\begin{array} {llll}
-k_{11}/k_{12}&-1 & e^{-i\gamma}/k_{12}&0\\
e^{i\gamma}/k_{12}&0&-k_{22}/k_{12}&1\end{array}\right]\in
\mathcal{O}^{\mathbb{C}}_{1,4}\cap\mathcal{O}^{\mathbb{C}}_{2,4}. $$
Since $k_{11}-\hat f_0k_{12}=0$,  $\mathbf{A}_1\in\mathcal{B}_{1,4}\cap\mathcal{B}_{2,4}$. Hence,
the conclusion holds by applying Theorems 4.1--4.2 to $\mathbf{A}_1$.
\medskip

Let $\lambda_n$ be restricted in $\mathcal{B}_{C} $ for each $ 0\leq n\leq N-1$. The following result is  a direct consequence of Corollary 4.3.
\medskip

\noindent{\bf Corollary 4.4.} {\it Fix a difference equation $\hat{\pmb\omega}$. Then for any given  $\mathbf{A}_1\in{\mathcal{B}_{C_1}}$, the $n$-th eigenvalue functions $\lambda_n,$ $0\leq n\leq N-1$, have the following asymptotic behaviors  near $\mathbf{A}_1$:
$$\begin{array} {llll}\lim\limits_{{\mathcal{B}_{C_1}^-}\ni\mathbf{A}\rightarrow\mathbf{A}_1}\lambda_0(\mathbf{A})=-\infty,\;
\lim\limits_{{\mathcal{B}_{C_1}^-}\ni\mathbf{A}\rightarrow\mathbf{A}_1}\lambda_n(\mathbf{A})=\lambda_{n-1}(\mathbf{A}_1), \;\; 1\leq n \leq N-1,\\[1.0ex]
\lim\limits_{{\mathcal{B}_{C_1}^+}\cup{\mathcal{B}_{C_1}}\ni\mathbf{A}\rightarrow\mathbf{A}_1}\lambda_n(\mathbf{A})=\lambda_{n}(\mathbf{A}_1), \;\; 0\leq n \leq N-2,\lim\limits_{{\mathcal{B}_{C_1}^+}\ni\mathbf{A}\rightarrow\mathbf{A}_1}\lambda_{N-1}(\mathbf{A})=+\infty.
\end{array}\vspace{-0.2cm}$$
}
\vspace{-0.2cm}
\bigskip

\noindent{\bf 5. Continuity and discontinuity of the $n$-th eigenvalue function in  the space of the  SLPs}\medskip

In this section, the continuous and discontinuous dependence of the $n$-th eigenvalue function  on the  SLP (1.1)--(1.2) is discussed. Its continuity and discontinuity sets in $\Omega_N^{\mathbb R,+} \times \mathcal B^{\mathbb C}$ are given and its asymptotic behaviors  near a discontinuity point are completely characterized.

Now, we introduce the following notations:
 \vspace{-0.2cm}$$\begin{array} {cccc} \mathcal{P}_{1,4}:=\{(\pmb\omega,\mathbf{A})\in\Omega_N^{\mathbb R,+} \times \mathcal O_{1,4}^{\mathbb C}:a_{12}=1/f_0\},
\mathcal{P}_{2,4}:=\left\{(\pmb\omega,\mathbf{A})\in\Omega_N^{\mathbb R,+} \times \mathcal O_{2,4}^{\mathbb C}: a_{11}=-f_0\right\},\\[2.0ex]
\mathcal{P}_{1,3}:=\left\{(\pmb\omega,\mathbf{A})\in\Omega_N^{\mathbb R,+} \times \mathcal O_{1,3}^{\mathbb C}: (a_{12}-1/f_0)b_{22}=|z|^2\right\},\\[2.0ex]
\mathcal{P}_{2,3}:=\left\{(\pmb\omega,\mathbf{A})\in\Omega_N^{\mathbb R,+} \times \mathcal O_{2,3}^{\mathbb C}:(a_{11}+f_0)b_{22}=|z|^2\right\},\\[2.0ex]
\mathcal{P}_{5}:=\left\{(\pmb\omega,\mathbf{A})\in\Omega_N^{\mathbb R,+} \times \mathcal B^{\mathbb C}:\mathbf{A}=\left [\begin{array} {cccc}1&1/ f_0&0&0\\
0&0&1&0\end{array}  \right ]\right\},\\[2.0ex]
\mathcal{P}_{1,3r}:=\left\{(\pmb\omega,\mathbf{A})\in\mathcal{P}_{1,3}: a_{12}-1/f_0\geq 0,b_{22}\geq0 \right\}\backslash\mathcal{P}_{5},
\\[2.0ex]
\mathcal{P}_{1,3l}:=\left\{(\pmb\omega,\mathbf{A})\in\mathcal{P}_{1,3}: a_{12}-1/f_0\leq 0,b_{22}\leq0 \right\}\backslash\mathcal{P}_{5},
\end{array}\vspace{-0.2cm}$$
\vspace{-0.2cm}$$\begin{array} {cccc}
\mathcal{P}^+_{1,3r}:=\left\{(\pmb\omega,\mathbf{A})\in\Omega_N^{\mathbb R,+} \times \mathcal O_{1,3}^{\mathbb C}: a_{12}-1/f_0\geq 0,b_{22}\geq0, (a_{12}-1/f_0)b_{22}>|z|^2\right\},\\[2.0ex]
\mathcal{P}^+_{1,3l}:=\left\{(\pmb\omega,\mathbf{A})\in\Omega_N^{\mathbb R,+} \times \mathcal O_{1,3}^{\mathbb C}: a_{12}-1/f_0\leq 0,b_{22}\leq0, (a_{12}-1/f_0)b_{22}>|z|^2\right\},\\[2.0ex]
\mathcal{P}_{2,3r}:=\left\{(\pmb\omega,\mathbf{A})\in\mathcal{P}_{2,3}:a_{11}+f_0\geq 0,b_{22}\geq 0\right\}\backslash\mathcal{P}_{5},\\[2.0ex]
\mathcal{P}_{2,3l}:=\left\{(\pmb\omega,\mathbf{A})\in\mathcal{P}_{2,3}:a_{11}+f_0\leq 0,b_{22}\leq 0\right\}\backslash\mathcal{P}_{5}, \end{array}$$
 $\mathcal{P}_{1,4}^+$ and $\mathcal{P}_{2,4}^+$, and $\mathcal{P}_{1,4}^-$, $\mathcal{P}_{2,4}^-$, $\mathcal{P}_{1,3}^-$, and $\mathcal{P}_{2,3}^-$ are defined similarly as $\mathcal{E}^+$ and
 $\mathcal{E}^-$, respectively;
$\mathcal{P}_{2,3r}^+$ and $\mathcal{P}_{2,3l}^+$ are defined similarly as
 $\mathcal{P}_{1,3r}^+$ and $\mathcal{P}_{1,3l}^+$, respectively.
 Then $\Omega_N^{\mathbb R,+} \times \mathcal O_{i,4}^{\mathbb C}=\mathcal{P}_{i,4}^+\cup\mathcal{P}_{i,4}\cup\mathcal{P}_{i,4}^-$, $\Omega_N^{\mathbb R,+} \times \mathcal O_{i,3}^{\mathbb C}=
\mathcal{P}_{i,3r}^+\cup\mathcal{P}_{i,3}\cup\mathcal{P}_{i,3}^-\cup\mathcal{P}_{i,3l}^+$, $\mathcal{P}_{i,3}=\mathcal{P}_{i,3r}\cup\mathcal{P}_5\cup\mathcal{P}_{i,3l}$ and $\mathcal{P}_{i,3r}\cap\mathcal{P}_{i,3l}=\varnothing$, where $i=1,2$.

Let $\mathcal{P}:=\cup_{i=1}^2\cup_{j=3}^4\mathcal{P}_{i,j}$. Note that  $\theta({\pmb\omega},\mathbf{A})=0$ if and only if $({\pmb\omega},\mathbf{A})\in\mathcal{P}$, where $\theta$ is defined by (2.5). Thus, $\mathcal{P}$ consists of  the SLPs  that have  less than $N$ eigenvalues. \vspace{0.1cm}
\medskip

\noindent{\bf Theorem 5.1. }{\it Each
$(\pmb\omega,
\mathbf{A})\in\mathcal{P}_{1,4}$ has exactly $N-1$ eigenvalues and  each $(\pmb\omega,
\mathbf{A})\in\left(\Omega_N^{\mathbb R,+} \times \mathcal O_{1,4}^{\mathbb C}\right)\backslash\mathcal{P}_{1,4}$ has exactly $N$ eigenvalues $\lambda_n(\pmb\omega,
\mathbf{A})$, $0\leq n \leq N-1$, which satisfy that
\begin{itemize}\vspace{-0.3cm}
\item[{\rm (i)}] $\lambda_n(\pmb\omega,
\mathbf{A})$ are continuous in $ \left(\Omega_N^{\mathbb R,+} \times \mathcal O_{1,4}^{\mathbb C}\right)\backslash\mathcal{P}_{1,4}$;\vspace{-0.3cm}
\item[{\rm (ii)}] $\lambda_n(\pmb\omega,
\mathbf{A})$ are not continuous at each point of $\mathcal{P}_{1,4}$ and have the following asymptotic behaviors  near any given $(\pmb\omega_0,
\mathbf{A}_0)\in \mathcal{P}_{1,4}$:
\vspace{-0.3cm}$$
\lim\limits_{\mathcal{P}_{1,4}^-\cup\mathcal{P}_{1,4}\ni(\pmb\omega,\mathbf{A})\rightarrow(\pmb\omega_0,\mathbf{A}_0)}\lambda_n(\pmb\omega,\mathbf{A})=\lambda_{n}(\pmb\omega_0,\mathbf{A}_0), \;\; 0\leq n \leq N-2,\eqno(5.1)\vspace{-0.7cm}$$
\vspace{-0.7cm}$$\vspace{-0.2cm}\lim\limits_{\mathcal{P}_{1,4}^-\ni(\pmb\omega,\mathbf{A})\rightarrow(\pmb\omega_0,\mathbf{A}_0)}\lambda_{N-1}(\pmb\omega,\mathbf{A})=+\infty,
\lim\limits_{\mathcal{P}_{1,4}^+\ni(\pmb\omega,\mathbf{A})\rightarrow(\pmb\omega_0,\mathbf{A}_0)}\lambda_0(\pmb\omega,\mathbf{A})=-\infty,\;\eqno(5.2)
\vspace{-0.5cm}$$
\vspace{-0.5cm}$$\lim\limits_{\mathcal{P}_{1,4}^+\ni(\pmb\omega,\mathbf{A})\rightarrow(\pmb\omega_0,\mathbf{A}_0)}\lambda_n(\pmb\omega,\mathbf{A})=\lambda_{n-1}(\pmb\omega_0,\mathbf{A}_0), \;\; 1\leq n \leq N-1,\eqno(5.3)\vspace{-0.5cm}$$
 and consequently, $\lambda_n$ restricted in $\mathcal{P}_{1,4}^-\cup\mathcal{P}_{1,4}$  is continuous  for each $0\leq n \leq N-2$.
\end{itemize}}

\noindent{\bf Proof.}
By  Lemma 2.4,   each $(\pmb\omega,\mathbf{A})\in\left(\Omega_N^{\mathbb R,+} \times\mathcal O_{1,4}^{\mathbb C}\right)\backslash\mathcal{P}_{1,4}$ has exactly $N$ eigenvalues, and each
$(\pmb\omega,\mathbf{A})\in\mathcal{P}_{1,4}$ has exactly $N-1$ eigenvalues.
Since
 $\left(\Omega_N^{\mathbb R,+} \times\mathcal O_{1,4}^{\mathbb C}\right)\backslash\mathcal{P}_{1,4}$ is an open subset of $\Omega_N^{\mathbb R,+} \times\mathcal B^{\mathbb C}$,
 by Theorem 2.1
$\lambda_n(\pmb\omega,
\mathbf{A})$  is continuous in $\left(\Omega_N^{\mathbb R,+} \times\mathcal O_{1,4}^{\mathbb C}\right)\backslash\mathcal{P}_{1,4}$
for each $0\leq n \leq N-1$.

The rest is to show that (5.1)--(5.3) hold for any given $(\pmb\omega_0,\mathbf{A}_0)\in\mathcal{P}_{1,4}$. Let $(r_1,r_2)$ be a finite interval such that $\lambda_j(\pmb\omega_0,\mathbf{A}_0)\in(r_1,r_2)$ for all $0\leq n\leq N-2$.
 Then by Lemma 2.5 there exists a  neighborhood $\mathcal{Y}$ of $(\pmb\omega_0,\mathbf{A}_0)$ in $\Omega_N^{\mathbb R,+} \times\mathcal O_{1,4}^{\mathbb C} $  such that each $(\pmb\omega,\mathbf{A})\in \mathcal{Y}$
has exactly $N-1$ eigenvalues in $[r_1,r_2]$ that  are all in $(r_1,r_2)$.
Let
\vspace{-0.2cm}$$\mathcal{Y}^-:= \mathcal{P}_{1,4}^-\cap \mathcal{Y},\;\mathcal{Y}^0:= \mathcal{P}_{1,4}\cap \mathcal{Y}, \;\mathcal{Y}^+:= \mathcal{P}_{1,4}^+\cap \mathcal{Y}.\vspace{-0.2cm}$$
 Note that
  $\mathcal{Y}$ can be chosen  such that $\mathcal{Y}^-$, $\mathcal{Y}^0$ and $\mathcal{Y}^+$ are connected.
 Since each $({\pmb\omega},\mathbf{A})\in\mathcal{Y}^+$ has exactly $N$ eigenvalues, it has exactly one
 eigenvalue, denoted by $\hat \lambda({\pmb\omega},\mathbf{A})$,  outside $[r_1,r_2]$.  By (i) of Theorem 2.3, either $L:=\{\hat \lambda({\pmb\omega},\mathbf{A}
):({\pmb\omega},\mathbf{A})\in\mathcal{Y}^+\}\subset(-\infty,r_1)$ or $L\subset(r_2, +\infty)$. For the  fix $\pmb\omega_0$,
there exists an $({\pmb\omega}_0,\mathbf{A}_1)\in\mathcal{Y}^+$ such that ${\lambda}_0({\pmb\omega}_0,\mathbf{A}_1)<r_1$ by the first relation in (4.2)  in Theorem 4.1. Hence,
 $\hat{\lambda}({\pmb\omega}_0,\mathbf{A}_1)={\lambda}_0({\pmb\omega}_0,\mathbf{A}_1)$.  Thus, again by (i) of Theorem 2.3,
 $L=\{\lambda_0({\pmb\omega},\mathbf{A}):({\pmb\omega},\mathbf{A})\in\mathcal{Y}^+\}\subset(-\infty,r_1)$.
Then, by  (ii) of Theorem 2.3, the second relation in  (5.2) holds.
With a similar argument, one can show that the first relation in  (5.2) holds.
  Note that $\mathcal{Y}^0$ is connected and  each
 $({\pmb\omega},\mathbf{A})\in\mathcal{Y}^0$ has exactly $N-1$ eigenvalues. So $\lambda_n$ restricted in $\mathcal{Y}^0$  is continuous and
  locally a continuous eigenvalue branch for each $0\leq n \leq N-2$  by Theorem 2.1. This, together with Theorem 2.2, yields that (5.1) and (5.3) hold. This completes the proof.
\medskip

With a similar method used in the proof of Theorem 5.1, one  can show that Theorems 5.2--5.4  hold.\medskip

\noindent{\bf Theorem 5.2. }{\it  Similar results in Theorem {\rm 5.1} hold for $\mathcal O_{1,4}^{\mathbb C}$, $\mathcal{P}_{1,4}$,  $\mathcal{P}_{1,4}^+$ and $\mathcal{P}_{1,4}^-$ replaced by
 $\mathcal O_{2,4}^{\mathbb C}$, $\mathcal{P}_{2,4}$, $\mathcal{P}_{2,4}^+$ and $\mathcal{P}_{2,4}^-$, separately.}
\medskip

\noindent{\bf Theorem 5.3. }{\it Each
$(\pmb\omega,
\mathbf{A})\in\mathcal{P}_5$ has exactly $N-2$ eigenvalues, and  each
$(\pmb\omega,
\mathbf{A})\in\mathcal{P}_{1,3}\backslash\mathcal{P}_5$ has exactly $N-1$ eigenvalues, and each $(\pmb\omega,
\mathbf{A})\in\left(\Omega_N^{\mathbb R,+} \times \mathcal O_{1,3}^{\mathbb C}\right)\backslash\mathcal{P}_{1,3}$ has exactly $N$ eigenvalues
 $\lambda_n(\pmb\omega,
\mathbf{A})$, $0\leq n \leq N-1$, which satisfy that
\begin{itemize}\vspace{-0.2cm}
\item[{\rm (i)}] $\lambda_n(\pmb\omega,
\mathbf{A})$  are continuous in $ \left(\Omega_N^{\mathbb R,+} \times \mathcal O_{1,3}^{\mathbb C}\right)\backslash\mathcal{P}_{1,3}$;\vspace{-0.3cm}
\item[{\rm (ii)}] $\lambda_n(\pmb\omega,
\mathbf{A})$ are not continuous at each point of $\mathcal{P}_{1,3}$ and  have the following asymptotic behaviors  near any given $(\pmb\omega_0,
\mathbf{A}_0)\in \mathcal{P}_{1,3r}$:
\vspace{-0.2cm}$$\begin{array} {llll}
\lim\limits_{\mathcal{P}_{1,3}^-\cup\mathcal{P}_{1,3r}\ni(\pmb\omega,\mathbf{A})\rightarrow(\pmb\omega_0,\mathbf{A}_0)}\lambda_{n}(\pmb\omega,\mathbf{A})=\lambda_{n}(\pmb\omega_0,\mathbf{A}_0), \;\; 0\leq n \leq N-2,\\[1.0ex]
\lim\limits_{\mathcal{P}_{1,3}^-\ni(\pmb\omega,\mathbf{A})\rightarrow(\pmb\omega_0,\mathbf{A}_0)}\lambda_{N-1}(\pmb\omega,\mathbf{A})=+\infty,
\lim\limits_{\mathcal{P}_{1,3r}^+\ni(\pmb\omega,\mathbf{A})\rightarrow(\pmb\omega_0,\mathbf{A}_0)}\lambda_{0}(\pmb\omega,\mathbf{A})=-\infty,\;\\[1.0ex]
\lim\limits_{\mathcal{P}_{1,3r}^+\ni(\pmb\omega,\mathbf{A})\rightarrow(\pmb\omega_0,\mathbf{A}_0)}\lambda_{n}(\pmb\omega,\mathbf{A})
=\lambda_{n-1}(\pmb\omega_0,\mathbf{A}_0),\;\;1\leq n \leq N-1;
\end{array}\vspace{-0.2cm}$$

near any given $(\pmb\omega_0,
\mathbf{A}_0)\in \mathcal{P}_{1,3l}$:
\vspace{-0.2cm}$$\begin{array}
{llll}\lim\limits_{\mathcal{P}_{1,3}^-\ni(\pmb\omega,\mathbf{A})\rightarrow(\pmb\omega_0,\mathbf{A}_0)}\lambda_{0}(\pmb\omega,\mathbf{A})=-\infty,\\[1.0ex]
\lim\limits_{\mathcal{P}_{1,3}^-\ni(\pmb\omega,\mathbf{A})\rightarrow(\pmb\omega_0,\mathbf{A}_0)}\lambda_{n}(\pmb\omega,\mathbf{A})=\lambda_{n-1}(\pmb\omega_0,\mathbf{A}_0), \;\; 1\leq n \leq N-1,\\[1.0ex]
\lim\limits_{\mathcal{P}_{1,3l}^+\cup\mathcal{P}_{1,3l}\ni(\pmb\omega,\mathbf{A})\rightarrow(\pmb\omega_0,\mathbf{A}_0)}\lambda_{n}(\pmb\omega,\mathbf{A})
=\lambda_{n}(\pmb\omega_0,\mathbf{A}_0),\;\;0\leq n \leq N-2,\\[1.0ex]
\lim\limits_{\mathcal{P}_{1,3l}^+\ni(\pmb\omega,\mathbf{A})\rightarrow(\pmb\omega_0,\mathbf{A}_0)}\lambda_{N-1}(\pmb\omega,\mathbf{A})=+\infty;
\end{array}\vspace{-0.2cm}$$
 and near any given $(\pmb\omega_0,
\mathbf{A}_0)\in\mathcal P_{5}$:
 \vspace{-0.2cm}$$\begin{array} {llll}\lim\limits_{\mathcal{P}_{1,3r}^+\cup\mathcal{P}_{1,3r}\cup\mathcal{P}_{1,3}^-\ni(\pmb\omega,\mathbf{A})\rightarrow(\pmb\omega_0,\mathbf{A}_0)}
 \lambda_{0}(\pmb\omega,\mathbf{A})=-\infty,
\lim\limits_{\mathcal{P}_{1,3r}^+\ni(\pmb\omega,\mathbf{A})\rightarrow(\pmb\omega_0,\mathbf{A}_0)}\lambda_{1}(\pmb\omega,\mathbf{A})=-\infty,\\[1.0ex]
\lim\limits_{\mathcal{P}_{1,3l}^+\cup\mathcal{P}_{1,3l}\cup\mathcal{P}_{5}\ni(\pmb\omega,\mathbf{A})\rightarrow(\pmb\omega_0,\mathbf{A}_0)}
\lambda_{n}(\pmb\omega,\mathbf{A})=\lambda_{n}(\pmb\omega_0,\mathbf{A}_0),\;0\leq n \leq N-3,\\[1.0ex]
\lim\limits_{\mathcal{P}_{1,3r}\cup\mathcal{P}_{1,3}^-\ni(\pmb\omega,\mathbf{A})\rightarrow(\pmb\omega_0,\mathbf{A}_0)}\lambda_{n}(\pmb\omega,\mathbf{A})
=\lambda_{n-1}(\pmb\omega_0,\mathbf{A}_0),\;1\leq n \leq N-2,\\[1.0ex]
\lim\limits_{\mathcal{P}_{1,3r}^+\ni(\pmb\omega,\mathbf{A})\rightarrow(\pmb\omega_0,\mathbf{A}_0)}\lambda_{n}(\pmb\omega,\mathbf{A})=
\lambda_{n-2}(\pmb\omega_0,\mathbf{A}_0), \;2\leq n \leq N-1,\\[1.0ex]
\lim\limits_{\mathcal{P}_{1,3l}^+\cup\mathcal{P}_{1,3l}\ni(\pmb\omega,\mathbf{A})\rightarrow(\pmb\omega_0,\mathbf{A}_0)}\lambda_{N-2}(\pmb\omega,\mathbf{A})=+\infty,
\lim\limits_{\mathcal{P}_{1,3l}^+\cup\mathcal{P}_{1,3}^-\ni(\pmb\omega,\mathbf{A})\rightarrow(\pmb\omega_0,\mathbf{A}_0)}\lambda_{N-1}(\pmb\omega,\mathbf{A})=+\infty.
\end{array}\vspace{-0.2cm}$$
\end{itemize}
And consequently, $\lambda_n$  restricted in  $\mathcal{P}_{1,3}^-\cup\mathcal{P}_{1,3r}$  and   $\mathcal{P}_{1,3l}^+\cup\mathcal{P}_{1,3l}\cup\mathcal{P}_{5}$ is continuous for each  $0\leq n \leq N-3$, and
$\lambda_{N-2}$  restricted in $\mathcal{P}_{1,3}^-\cup\mathcal{P}_{1,3r}$ and $\mathcal{P}_{1,3l}^+\cup\mathcal{P}_{1,3l}$ is continuous.}\medskip

\noindent{\bf Theorem 5.4. }{\it Similar results in Theorem {\rm 5.3} hold for  $\mathcal{O}_{1,3}^{\mathbb{C}}$,  $\mathcal{P}_{1,3}$, $\mathcal{P}_{1,3r}$, $\mathcal{P}_{1,3l}$, $\mathcal{P}_{1,3}^-,$ $\mathcal{P}_{1,3r}^+$, and $\mathcal{P}_{1,3l}^+$ replaced by $\mathcal{O}_{2,3}^{\mathbb{C}}$, $\mathcal{P}_{2,3}$,
 $\mathcal{P}_{2,3r}$, $\mathcal{P}_{2,3l}$, $\mathcal{P}_{2,3}^-,$ $\mathcal{P}_{2,3r}^+$, and $\mathcal{P}_{2,3l}^+$, separately.}\medskip

Combining Theorems 5.1--5.4 yields the continuity and discontinuity sets of the $n$-th eigenvalue function in $\Omega_N^{\mathbb R,+} \times \mathcal B^{\mathbb C}$:\medskip

\noindent{\bf Theorem 5.5.} {\it  Each $(\pmb\omega,\mathbf{A})\in \mathcal{P}_5$
has exactly $N-2$ eigenvalues, and  each   $(\pmb\omega,\mathbf{A})\in\mathcal{P}\backslash\mathcal{P}_5$ has exactly $N-1$ eigenvalues, and  each $(\pmb\omega,\mathbf{A})\in\left(\Omega_N^{\mathbb R,+} \times \mathcal B^{\mathbb C}\right)\backslash\mathcal{P}$
 has exactly $N$ eigenvalues, and $\lambda_n$, $0\leq n\leq N-1$, are continuous in $ \left(\Omega_N^{\mathbb R,+} \times \mathcal B^{\mathbb C}\right)\backslash\mathcal{P}$
 and  not continuous at each point of $\mathcal{P}$. }\medskip

\noindent{\bf Remark 5.1.} If the  $n$-th eigenvalue function $\lambda_n$ is  simple at some $( \pmb\omega_0,\mathbf{A}_0)$ in the continuity set $\left(\Omega_N^{\mathbb R,+} \times \mathcal B^{\mathbb C}\right)\backslash\mathcal{P}$ for some $0\leq n\leq N-1$, the  results about differentiability of continuous eigenvalue branches in  [22, Theorems 4.2-4.4 and 4.7] can be applied to $\lambda_n$ in a neighborhood of $( \pmb\omega_0,\mathbf{A}_0)$.

\bigskip \noindent{\bf \large References}
\def\hang{\hangindent\parindent}
\def\textindent#1{\indent\llap{#1\enspace}\ignorespaces}
\def\re{\par\hang\textindent}
\noindent \vskip 3mm

\re{[1]} F. V. Atkinson,
Discrete and Continuous Boundary Problems, Academic
Press, New York, 1964.

\re{[2]} M. Bohner, O. Dosly, The discrete Pr\"{u}fer transformation, Proc. Amer. Math. Soc. 129 (2001) 2715--2725.

\re{[3]} M. Bohner, O. Dosly, W. Kratz, Sturmian and spectral theory for discrete symplectic
systems, Trans. Amer. Math. Soc. 361 (2009)
3109--3123.

\re{[4]} S. Clark, A spectral analysis for self-adjoint operators generated
by a class of second order difference equations, J. Math.
Anal. Appl. 197 (1996) 267--285.

\re{[5]} S. Clark, F. Gesztesy, W. Renger, Trace formulas and Borg-type
theorems for matrix-valued Jacobi and Dirac finite difference
operators, J. Differ. Equ. 219 (2005)
144--182.

\re{[6]} R. Courant, D. Hilbert, Methods of Mathematical Physics, Interscience Publishers, New York, 1953.

\re{[7]} M. Eastham, Q. Kong, H. Wu, A. Zettl,  Inequalities among eigenvalues of Sturm-Liouville problems, J. Inequal. Appl. 3 (1999) 25--43.

\re{[8]} W. N. Everitt, M. M\"{o}ller, A. Zettl, Discontinuous dependence of the $n$-th Sturm-Liouville eigenvalue, Int. Ser. Numer. Math. 123 (1997) 145--150.

\re{[9]} A. Jirari,  Second-order Sturm-Liouville difference equations and
orthogonal polynomials, Mem. Amer. Math. Soc. 113 (1995).

\re{[10]} Q. Kong, H. Wu, A. Zettl, Dependence of the $n$-th
Sturm-Liouville eigenvalue on the problem, J. Differ.
Equ. 156 (1999) 328--354.

\re{[11]} Q. Kong, H. Wu, A. Zettl, Geometric aspects of Sturm-Liouville
problems, I. Structures on spaces of boundary conditions,
Proc. Roy. Soc. Edinb. Sect. A Math. 130 (2000) 561--589.

\re{[12]} H. Lv, Y. Shi, Error estimate of eigenvalues of perturbed second-order discrete Sturm-
Liouville problems, Linear Algebra Appl. 430 (2009) 2389--2415.

\re{[13]} W. Peng, M. Racovitan, H. Wu, Geometric aspects
of Sturm-Liouville problems, V. Natural loops of boundary conditions
for monotonicity of eigenvalues and their applications, Pac. J. Appl. Math. 4 (2006) 253--273.

\re{[14]} J. Poeschel, E. Trubowitz, Inverse Spectral Theory, Academic
Press, New York, 1987.

\re{[15]} Y. Shi, Spectral theory for a class of discrete linear Hamiltonian
systems, J. Math. Anal. Appl. 289 (2004) 554--570.

\re{[16]} Y. Shi,  Weyl-Titchmarsh theory for a class of discrete linear
Hamiltonian systems, Linear Algebra Appl. 416 (2006)
452--519.

\re{[17]} Y. Shi, S. Chen,  Spectral theory of second--order vector
difference equations, J. Math. Anal. Appl. 239
(1999) 195--212.

\re{[18]} Y. Shi, S. Chen,  Spectral theory of higher-order discrete
vector Sturm-Liouville problems, Linear Algebra Appl.
323 (2001) 7--36.

\re{[19]} H. Sun, Y. Shi, Eigenvalues of second-order difference
equations with coupled boundary conditions, Linear Algebra
Appl. 414 (2006) 361--372.

\re{[20]} Y. Wang, Y. Shi,  Eigenvalues of second-order difference
equations with periodic and anti-periodic boundary conditions,
J. Math. Anal. Appl. 309 (2005) 56--69.

\re{[21]} A. Zettl, Sturm-Liouville Theory, Mathematical Surveys Monographs, vol. 121, Amer. Math. Soc., 2005.

\re{[22]} H. Zhu, S. Sun, Y. Shi, H. Wu, Dependence of discrete Sturm-Liouville eigenvalues on problems, submitted for publication.

\end{document}